\documentclass{gtart_h}

\def\ifplaintex{\expandafter\ifx\csname documentclass\endcsname\relax}

\ifplaintex 
\else
\fi

\def\gt{{\mathsurround=0pt\it $\cal G\mskip-2mu$eometry \&\ 
$\cal T\!\!$opology}}        

\def\gtm{{\mathsurround=0pt\it $\cal G\mskip-2mu$eometry \&\ 
$\cal T\!\!$opology $\cal M\mskip-1mu$onographs}}    

\def\gtp{{\mathsurround=0pt\it $\cal G\mskip-2mu$eometry \&\ 
$\cal T\!\!$opology $\cal P\!$ublications}}  

\def\recd{{\small Received:\qua\receiveddate\ifx\reviseddate\relax
\else\qquad Revised:\qua\reviseddate\fi\par}} 

\def\volumenumber#1{\def\thevolumenumber{#1}}
\def\volumeyear#1{\def\thevolumeyear{#1}}
\def\volumename#1{\def\thevolumename{#1}}
\def\papernumber#1{\def\thepapernumber{#1}}
\def\pagenumbers#1#2{\def\startpage{#1}\def\finishpage{#2}}
\def\published#1{\def\publishdate{#1}}
\def\received#1{\def\receiveddate{#1}}
\def\revised#1{\def\reviseddate{#1}}
\def\accepted#1{\def\accepteddate{#1}}

\def\asciiemail#1{\def\theasciiemail{#1}}

\long\def\asciiabstract#1{\long\def\theasciiabstract{#1}}
\def\asciikeywords#1{\def\theasciikeywords{#1}}

\let\\\par
\let\thevolumenumber\relax\let\thepapernumber\relax
\let\thevolumeyear\relax\let\startpage\relax
\let\finishpage\relax\let\publishdate\relax\let\receiveddate\relax
\let\reviseddate\relax\let\accepteddate\relax\let\theasciititle\relax
\let\theasciiauthors\relax
\let\theasciiabstract\relax\let\theasciikeywords\relax

\let\theerratum\relax\let\theasciiemail\relax
\let\theshortauthors\relax\let\theshorttitle\relax

\def\startpage{1}\def\finishpage{15}\def\thepapernumber{77}

\volumenumber{2}
\volumename{Proceedings of the Kirbyfest}
\volumeyear{1999}

\long\def\maketitlep{   

\count0=\startpage

\gtm\nl        
{\small Volume \thevolumenumber: \thevolumename\nl 
\ifx\theerratum\relax\else Erratum \erratumnumber\nl\fi
Pages \startpage--\finishpage\nl}

\vglue 0.1truein   

{\parskip=0pt\leftskip 0pt plus 1fil\def\\{\par\smallskip}{\ifplaintex\large
\else\Large\fi\bf\thetitle}\par\medskip}   
\vglue 0.05truein 

{\parskip=0pt\leftskip 0pt plus 1fil\def\\{\par}{\sc\theauthors}
\par\medskip}%
 
\vglue 0.03truein 

{\small\leftskip 25pt\rightskip 25pt{\bf Abstract}\stdspace\theabstract

{\bf AMS Classification}\stdspace\theprimaryclass
\ifx\thesecondaryclass\relax\else; \thesecondaryclass\fi\par
{\bf Keywords}\stdspace \thekeywords\par}\vglue 7pt

}   

\font\phead=cmsl9 scaled 950
\font=cmsl9 scaled 1050
\font\pnum=cmbx10 scaled 913
\font\lnum=cmbx10 
\font\pfoot=cmsl9 scaled 950
\font=cmsl9 scaled 1050
\ifplaintex
\headline{\vbox to 0pt{\vskip -4.5mm\line{\small\phead\ifnum
\count0=\startpage ISSN 1464-8997 (on line)
1464-8989 (printed) \hfill {\pnum\folio}\else\ifodd\count0\def\\{ }%
\ifx\theshorttitle\relax\thetitle\else\theshorttitle\fi\hfill{\pnum\folio}
\else\def\\{ and }{\pnum\folio}\hfill\ifx\theshortauthors\relax\theauthors
\else\theshortauthors\fi\fi\fi}\vss}}
\footline{\vbox to 0pt{\vglue 0mm\line{\small\pfoot\ifnum\count0=\startpage
Published \publishdate:\qua\copyright\ \gtp\hfill\else
\gtm, Volume \thevolumenumber\ (\thevolumeyear)\hfill\fi}\vss
}}
\else
\makeatletter
\def\@oddhead{{\small\lhead\ifnum\count0=\startpage ISSN 1464-8997 (on line)
1464-8989 (printed) \hfill {\lnum\number\count0}\else\ifodd\count0
\def\\{ }\ifx\theshorttitle\relax \thetitle \else\theshorttitle\fi\hfill
{\lnum\number\count0}\else\def\\{ and }{\lnum\number\count0}
\hfill\ifx\theshortauthors\relax 
\theauthors\else\theshortauthors\fi\fi\fi}}\def\@evenhead{@oddhead}
\def\@oddfoot{\small\lfoot\ifnum\count0=\startpage Published \publishdate:\qua\copyright\ \gtp\hfill\else
\gtm, Volume \thevolumenumber\ (\thevolumeyear)\hfill\fi}
\def\@evenfoot{@oddfoot}
\makeatother
\fi

\let\maketitlepage\maketitlep
\let\makeshorttitle\maketitlepage
\let\maketitle\maketitlepage

\newwrite\gtoutfile
\long\gdef\makeheadfile{  
{\def\\{, }\def\s{ }
\immediate\openout\gtoutfile head.xxx
\immediate\write\gtoutfile{Proxy-for: \ifx\theasciiauthors\relax
\theauthors\else\theasciiauthors\fi\s<\ifx\theasciiemail\relax\theemail\else\theasciiemail\fi>}
\immediate\write\gtoutfile{\noexpand\\}
\immediate\write\gtoutfile{Authors: \ifx\theasciiauthors\relax
\theauthors\else\theasciiauthors\fi}
{\def\\{ }\immediate\write\gtoutfile{Title: \ifx\theasciititle\relax
\thetitle\else\theasciititle\fi}}
\immediate\write\gtoutfile{Subj-class: GT or SG, GR etc}
\immediate\write\gtoutfile{MSC-class: \theprimaryclass\ifx\thesecondaryclass\relax\else, \thesecondaryclass\fi}
\immediate\write\gtoutfile{Journal-ref: Geom. Topol. Monogr. \thevolumenumber\s
(\thevolumeyear) \startpage-\finishpage}
\immediate\write\gtoutfile{Comments: Published by Geometry and Topology Monographs at}
\immediate\write\gtoutfile{\s\s\s  http://www.maths.warwick.ac.uk/gt/GTMon\thevolumenumber/paper\thepapernumber.abs.html}
\immediate\write\gtoutfile{\noexpand\\}
\immediate\write\gtoutfile{}
\ifx\theasciiabstract\relax
\immediate\write\gtoutfile{\theabstract}\else
\immediate\write\gtoutfile{\theasciiabstract}\fi
\immediate\write\gtoutfile{}
\immediate\write\gtoutfile{\noexpand\\}
\immediate\write\gtoutfile{}
\immediate\closeout\gtoutfile}}  

\def\maketitlepage{\maketitlep\makeheadfile}
\let\makeshorttitle\maketitlepage
\let\maketitle\maketitlepage

\volumenumber{7}
\volumename{Proceedings of the Casson Fest}
\volumeyear{2004}
\papernumber{4}
\pagenumbers{101}{134}
\received{4 December 2003}
\revised{24 July 2004}
\accepted{17 June 2004}
\published{18 September 2004}

\usepackage{amsmath,amssymb,graphicx}

\newtheorem{thm}{Theorem}

\newtheorem{lem}[thm]{Lemma}
\newtheorem{prop}[thm]{Proposition}
\newtheorem{conj}[thm]{Conjecture}
\theoremstyle{definition}
\newtheorem{defn}[thm]{Definition}
\newtheorem{rem}[thm]{Remark}

\renewcommand{\int}{\operatorname{int}}

\newcommand{\Z}{\mathbb{Z}}

\newcommand{\Q}{\mathbb{Q}}

\newcommand{\cT}{\mathcal{T}}

\newcommand{\cW}{\mathcal{W}}

\newcommand{\e}{\epsilon}

\def\theenumi{\roman{enumi}}

\begin{document}

\title{Whitney towers and the Kontsevich integral}
\author{Rob Schneiderman\\Peter Teichner}
\address{Courant Institute of Mathematical Sciences,
New York University\\251 Mercer Street, New York, NY 10012-1185, USA}
\address{Department of Mathematics, University of California\\Berkeley,
CA 94720-3840, USA}
\gtemail{\mailto{schneiderman@courant.nyu.edu},
\mailto{teichner@math.berkeley.edu}}
\asciiemail{schneiderman@courant.nyu.edu, teichner@math.berkeley.edu}
\begin{abstract}
We continue to develop an obstruction theory for embedding $2$--spheres
into $4$--manifolds in terms of Whitney towers. The proposed intersection
invariants take values in certain graded abelian groups generated
by labelled trivalent trees, and with relations well known from the
3--dimensional theory of finite type invariants.  Surprisingly, the
same exact relations arise in 4 dimensions, for example the Jacobi (or
IHX) relation comes in our context from the freedom of choosing Whitney
arcs. We use the finite type theory to show that our invariants agree
with the (leading term of the tree part of the) Kontsevich integral
in the case where the $4$--manifold is obtained from the $4$--ball by
attaching handles along a link in the $3$--sphere.
\end{abstract}

\asciiabstract{We continue to develop an obstruction theory for embedding 2-spheres
into 4-manifolds in terms of Whitney towers. The proposed intersection
invariants take values in certain graded abelian groups generated
by labelled trivalent trees, and with relations well known from the
3-dimensional theory of finite type invariants.  Surprisingly, the
same exact relations arise in 4 dimensions, for example the Jacobi (or
IHX) relation comes in our context from the freedom of choosing Whitney
arcs. We use the finite type theory to show that our invariants agree
with the (leading term of the tree part of the) Kontsevich integral
in the case where the 4-manifold is obtained from the 4-ball by
attaching handles along a link in the 3-sphere.}

\asciikeywords{2-sphere, 4-manifold, link concordance,
Kontsevich integral, Milnor invariants, Whitney tower}

\maketitlepage

\cl {\small\it Dedicated to Andrew Casson on the occasion of his 60th
birthday}

\section{Introduction}
Two of Andrew Casson's wonderful contributions to topology were
his work on {\em flexible handles} (now called {\em Casson
towers}) in $4$--manifolds, and his invariant for homology
$3$--spheres, counting representations into $SU(2)$. In this paper
we will describe an obstruction theory for disjointly embedding
collections of 2--spheres (or 2--disks with fixed boundary) into a
4--manifold that provides a connection between these two aspects
of Casson's work. This connection is somewhat indirect, otherwise
our paper would be called {\em Casson towers and the Casson
invariant}. In other words, we shall switch from Casson towers to
Whitney towers, and from the Casson invariant to the Kontsevich
integral. It would be very satisfying to find a more
straightforward relationship between Casson's two contributions.

To explain the connection, recall that the Casson invariant is the
lowest order (nontrivial) {\em finite type invariant} of homology
$3$--spheres. These finite type invariants take values in certain
graded abelian groups generated by trivalent graphs. Being {\em
invariants}, they measure the {\em uniqueness} of $3$--manifolds
or links in $3$--manifolds. We shall explain how similar graphs,
better, unitrivalent trees, arise in {\em existence} questions for
$4$--manifolds or surfaces in $4$--manifolds. It is not totally
surprising that raising the dimension by one takes uniqueness to
existence questions, after all an isotopy of, say, a knot in a
$3$--manifold $M^3$ is nothing but a certain annulus in the
$4$--manifold $M \times I$. However, the details of such a
translation from one dimension to the next are not at all obvious.

In the easiest setting one would like to find obstructions for
making the images of maps $A_i\co(D^2,S^1)\to (X^4,\partial X)$ {\em
disjoint}, without changing the homotopy classes (and without
trying to embed the $A_i$). In fact, Casson's main Embedding
Theorem in \cite[Lecture 1]{Ca} is an example of a special case of
this problem: Casson showed that if $X$ is simply connected, all
intersection numbers between the $A_i$ vanish, and $A_i$ have {\em
algebraic dual spheres}, then the problem has a positive solution.
He used inverses of the Whitney move, now known as Casson or
finger moves, to introduce many self-intersections, while
trivializing the fundamental group of the complement of one disk
at a time (and hence enabling the other disks to be mapped {\em
disjointly}). He then went on to construct Casson towers (with
prescribed boundary circles) by iterating the procedure
indefinitely, using the fact that the complement of a finite
height Casson tower can be made simply connected. These ideas
inspired Mike Freedman who proved in \cite{F} that a neighborhood
of a Casson tower actually contains an embedded flat disk.

The presence of algebraic dual spheres in Casson's theorem comes
from the fact that the proposed application was to the s-cobordism
theorem and to the exactness of the surgery sequence in
dimension~4. Indeed,  Freedman's theorem implies these results in
the topological category (for {\em good} fundamental groups).

There is a more general context in which disjoint maps of disks or
spheres can be constructed, namely in the presence of a {\em
non-repeated Whitney tower} (of sufficiently high order), see
Theorem~\ref{thm:nonrepeating} below and \cite{ST2}. The first
order stage of this Whitney tower is guaranteed by the vanishing
of the intersection numbers whereas the existence of the higher
order stages are obstructed by our new proposed invariants. They
take values in certain graded abelian groups generated by
trivalent trees, which are basically the spines of the Whitney
towers. The difference between a Casson tower and a Whitney tower
is that in the latter, fewer disks are attached at each stage: In
a Casson tower, every intersection point $p$ leads to a new disk
(with boundary an arc leaving on one sheet at $p$ and arriving at
the other sheet), whereas a Whitney tower only has a new disk for
certain {\em pairs} of intersection points. In particular, it is
usually only possibly to find Casson towers in simply connected
$4$--manifolds, whereas Whitney towers are not restricted by the
fundamental group. In fact, in our theory the fundamental group
leads to a decoration of the trivalent trees in question, thus
giving a much bigger variety of possible obstructions. In
addition, Freedman's reimbedding theorem shows that a Casson tower
of height 3 already contains an embedded flat disk. However, there
are Whitney towers of arbitrary order {\em not} containing disks,
which explains the use of these ``weaker'' towers in an
obstruction theory.

Our Theorem~\ref{thm:nonrepeating} implies Casson's result because
algebraic dual spheres can be used to construct non-repeating
Whitney towers of arbitrary order. This is already implicit in
\cite{FQ}, so our main contribution is a theory {\em in the
absence} of algebraic dual spheres. For example, this applies to
concordance questions for links in $3$--space. In this context we
prove in Theorem~\ref{thm:Milnor} below that our invariants agree
rationally with (the leading term of) the tree part of the
Kontsevich integral, which is the universal finite type
concordance invariant \cite{HM}. This relates our obstruction
theory to the finite type theory and, in particular, to the Casson
invariant. It should be mentioned here that Habegger and Masbaum
show in \cite{HM} that (the leading term of) the tree part of the
Kontsevich integral carries exactly the same information as
Milnor's $\overline{\mu}$--invariants which were first observed to
be concordance invariants by Casson in \cite{Ca2}. Reversing the
logic, we have found a $4$--dimensional geometric interpretation
of this part of the Kontsevich integral, in terms of higher order
intersections among Whitney disks. See  \cite{CT2} for an
interpretation in terms of {\em gropes} in $3$--dimensions which
is stronger in the sense that it works for (the leading term of)
the Kontsevich integral, not just of the tree part.

At the time of writing, the setting of Theorem~\ref{thm:Milnor} is
actually the only case where we have a proof that our intersection
invariant is independent of the choice of a Whitney tower, but see
Conjecture~\ref{conj:well defined}. What we do prove in
Theorem~\ref{thm:build-tower} is that the vanishing of our
intersection invariant for a Whitney tower of order~$n$ enables
one to build a Whitney tower of the next order~$(n+1)$. In that
sense, we are producing an obstruction theory since disjointly
embedded sheets $A_i$ allow Whitney towers of arbitrary order.

We close this introduction by pointing out that the Whitney towers
used in this paper are generalizations of the ones in \cite{COT}
in that disks of higher order are here allowed to intersect
previous stages, as long as these intersection points are paired
up by Whitney disks (up to the desired order). In our language, the distinction is made in
terms of saying that these Whitney towers have an {\em order}
whereas the Whitney towers of \cite{COT} (where different order
Whitney disks don't intersect) have a {\em height}. This is the
precise analogue of {\em class} versus {\em height} in the theory
of gropes, see eg \cite{T}, ultimately coming from the
distinction between the lower central series and the derived
series of a group. The latter explains why Whitney towers with
a height carry more subtle information. In fact, they are {\em
not} related to the usual finite type theory and hence it is much
more difficult to define an obstruction theory. At present, such a
theory only exists for knot concordance \cite{COT}, \cite{CT} (using von
Neumann signatures to prove nontriviality) and it would be
extremely interesting to develop it more generally, ie in the
context of $2$--spheres in $4$--manifolds.

\section{Statement of results}
We continue to develop the obstruction theory for embedding $2$--spheres
into $4$--manifolds started in \cite{ST}. To fix notation, let $X$ be a
$4$--manifold and $A_1,\dots, A_m$ be generic immersions of $2$--spheres
(or $2$--disks with fixed boundary) into $X$. We shall work in the smooth
setting, even though the techniques of \cite{FQ} allow a generalization
of our work to locally flat surfaces in a topological manifold. The goal
is to construct obstructions for changing the $A_i$, in their regular
homotopy class, to embeddings with disjoint images. This is already a
very interesting problem for $m=1$ but we shall not restrict to this case.

The first, well known, invariants are the Wall intersection ``numbers''
\cite{W}
$$\lambda(A_i,A_j) = \sum_{p\in A_i\cap A_j} \e_p\cdot g_p \quad
  \in \Z\pi, \quad \pi:=\pi_1X.$$
These count how often $A_i$ and $A_j$ intersect algebraically, including
a group element $g_p\in\pi$ and a sign $\e_p$ for each intersection
point. Similarly, there are self-intersection numbers $\mu(A_i)$ which
are well defined only in a certain quotient of the group ring, see
below. Recall that in higher dimensions (where $A_i$ are $k$--spheres,
$k>2$ and $X$ is $2k$--dimensional) the vanishing of these invariants
implies that after a finite sequence of Whitney moves \cite{Wh} the $A_i$
can be represented by disjoint embeddings.  In dimension 4, there are well
known problems to this procedure  (since $2+k=2k$ for $k=2$), the most
important one being that, generically, the Whitney disks intersect the
$2$--spheres $A_i$. The first precise statement concerning the failure
of the Whitney trick in dimension~4 was given by Kervaire and Milnor
in \cite{KM}.

In \cite{ST} we assumed that these primary  intersection numbers
vanish which means geometrically that all intersections and
self-intersections can be paired by Whitney disks: For each pair
of intersection points between $A_i$ and $A_j$ (if $i=j$ these are
self-intersections), choose one {\em Whitney arc} on $A_i$ and one
on $A_j$ connecting these two points. Since the fundamental group
is controlled in Wall's invariant, the two Whitney arcs together
form a null homotopic circle in the ambient $4$--manifold, which
hence bounds a disk, the {\em Whitney disk}. Using a choice for
such disks, one for each pair of intersection points, we
constructed a secondary invariant
$$\tau(A_i,A_j,A_k) \in \Z\pi \times \Z\pi /\dots$$
which measures how the Whitney disks intersect the spheres $A_i$. Here the
indices $i,j,k$ may be repeated, obtaining several slightly distinct
geometrical cases just like for Wall's invariants. We recall that
by standard procedures the Whitney disks can always be assumed to be
disjointly embedded (and framed), and that the only thing which hinders a
successful Whitney move is the fact that they are in general not disjoint
from the original spheres $A_i$.

We will first explain a way to unify the above invariants, then
suggest a vast generalization and finally discuss a relation to
Milnor invariants and the Kontsevich integral (for classical links).
For this purpose, assume that the $A_i$ intersect and self-intersect
generically, and call the collection $A_1,\dots,A_m$ a {\em Whitney
tower of order~0}. Similarly, if Wall's invariants vanish, and one has
chosen generic Whitney disks $W_I$ which pair all intersections and
self-intersections of the $A_i$ then one obtains a {Whitney tower of
order~1}. If the $\tau$--invariants vanish, then one can chose Whitney
disks for all the intersections of the $A_i$ with the $W_I$ to obtain
a {Whitney tower of order~2}. This procedure can be continued and
we give a precise definition of a {\em Whitney tower of order~n} in
Section~\ref{W-tower-tau-sec}. This definition includes orientations of
all the surfaces $A_i,W_I,\dots$ in the tower, as well as base points
on these surfaces together with whiskers connecting these base points
to the base point of $X$.

\begin{figure}[ht!]
\centerline{\includegraphics[scale=.4]{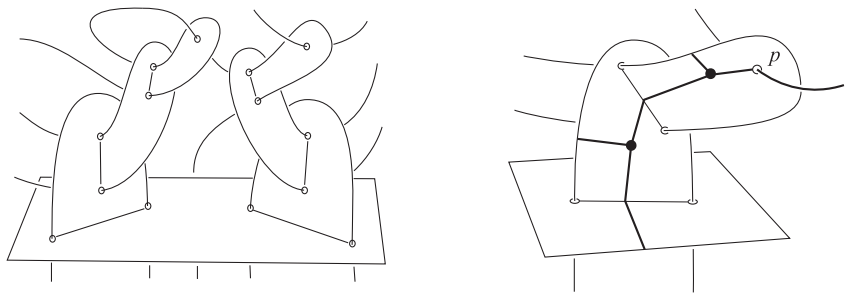}}
\caption{Part of a Whitney tower (left), and part of the
  unitrivalent tree $t_p$ associated to an unpaired intersection point
  $p$ in a Whitney tower (right).}
\label{W-tower5-int-point-tree-fig}
\end{figure}

\subsection{The intersection tree $\tau_n(\cW)$}
\label{tau-intro}

Our first observation is that one can canonically associate to
each {\em unpaired} intersection point $p$ of a Whitney tower
$\cW$ a {\em decorated}  unitrivalent tree $t_p$ of {\em
order}~$n$. The order is  the number of trivalent vertices and the
decoration is as follows: the univalent vertices of $t_p$ are
labelled by the $A_i$ or more abstractly, by $i\in\{1,\dots,m\}$,
the edges are labelled by elements from the fundamental group
$\pi$, and the edges and trivalent vertices are oriented. The tree
$t_p$ sits naturally as a subset of $\cW$
(Figure~\ref{W-tower5-int-point-tree-fig}, details in
Section~\ref{W-tower-tau-sec}) with each trivalent vertex lying in
a Whitney disk and each univalent vertex lying in some $A_i$. Each
edge of $t_p$ is a sheet-changing path between vertices in
adjacent surfaces, with the group element labelling the edge
determined by the loop formed from the path together with the
whiskers on the adjacent surfaces. For example, in a Whitney tower
of order~0, any intersection point $p$ between $A_i$ and $A_j$ has
order~0 and gives a tree $t_p$ consisting of a single edge whose
univalent vertices (labelled by $i$ and $j$) correspond to
basepoints in $A_i$ and $A_j$. This edge is labelled by the group
element $g_p$ determined by a loop formed from the whiskers on
$A_i$ and $A_j$ together with a path that changes sheets at $p$
where the orientation of the edge corresponds to the direction of
the path. For intersection points of order~1 in an order--1 Whitney
tower, one gets decorated Y--trees with one trivalent vertex and
three univalent vertices labelled by $i,j,k$ (which can repeat).

The central point of this paper is that in an order--$n$ Whitney
tower $\cW$ the trees that correspond to the (unpaired) order--$n$
intersection points of $\cW$ represent a ``higher order''
obstruction to homotoping (rel boundary) the $A_i$ to disjoint
embeddings. Just like the intersection number $\lambda(A_i,A_j)$
is a sum over all intersection points between $A_i$ and $A_j$, we
define the {\em intersection tree} $\tau_n(\cW)$ of an order--$n$
Whitney tower $\cW$ to be
$$\tau_n(\cW):= \sum_p \e_p\cdot t_p \quad \in \cT_n(\pi,m).$$
The sum is taken over all order--$n$ intersection points $p$ in
$\cW$ and we consider this sum as taking values in the free
abelian group generated by (isomorphism classes of) decorated
trees as above, modulo several relations that are motivated
geometrically (explained briefly below and in detail in
Section~\ref{W-tower-tau-sec}, particularly
Section~\ref{subsec:T}). We denote this quotient by
$$\cT(\pi,m)=\bigoplus_{n=0}^\infty \cT_n(\pi,m),$$
where the order~$n$ is the number of trivalent vertices and the univalent
labels come from $\{1,\dots,m\}$, possibly repeated. If this index set
is undetermined (or unimportant) we shall just write $\cT_n(\pi)$.

The order--0 trees are just single edges and it turns out that
$$\cT_0(\pi,1) \cong \Z\pi/ \langle \bar g- g \rangle, \quad
\bar g:=w_1(g)\cdot g^{-1},$$
where $w_1\co\pi\to \Z/2$ is the first Stiefel--Whitney class of the
ambient 4--manifold. The quotient comes from the fact that an edge
with two identical labels has an additional symmetry which changes
the orientation of the edge. Moreover, our invariant $\tau_0$
gives exactly Wall's self-intersection invariant $\mu$. To get
Wall's intersection number $ \lambda(A_1,A_2)$ we just need to
evaluate $\tau_0$ in order~0 with {\em exactly} two labels~$1,2$.
The invariants $\tau$ from \cite{ST} are exactly $\tau_1(\cW)$ in
the various versions of $\cT_1(\pi)$, depending on the allowed
labels.

A short discussion of the relations in $\cT(\pi)$ is in order.
They reflect the various choices made in the construction of the
Whitney tower, as will be discussed in
Section~\ref{W-tower-tau-sec} (see also Figure~\ref{Relations-fig}
in Section~\ref{W-tower-tau-sec}). As a consequence, working {\em
modulo} these relations makes our intersection tree $\tau_n$ {\em
independent} of the choices below.
\begin{itemize}
\item Changing orientations on Whitney disks gives AS,
{\em antisymmetry} relations; they introduce a sign when the cyclic
ordering of a trivalent vertex is switched.
\item Changing the orientation of an edge changes the label $g$ to $\bar
g$, the OR {\em orientation} relation.
\item Changing the whiskers gives HOL, {\em holonomy} relations; they
multiply the labels of 3 edges coming into a
trivalent vertex by a group element.
\item Changing the choice of Whitney arcs, ie of the boundaries of
Whitney disks, gives the IHX relations.
\end{itemize}

The last type of relations, well known in dimension~3, is maybe the
most surprising aspect of our $4$--dimensional theory. We feel that
our explanation in terms of the indeterminacy of Whitney arcs is very
satisfying \cite{CST}. It should be pointed out that graded abelian
groups like $\cT(\pi)$ arose independently in the $3$--dimensional work
of Garoufalidis, Kricker and Levine \cite{GK}, \cite{GL}. They study
trivalent graphs (instead of unitrivalent trees) and $\pi$ is usually
a $3$--manifold group. In some form, the Kontsevich integral gives
invariants of links (or $3$--manifolds) with values in such graded
abelian groups. So these are invariants for the {\em uniqueness} of
$3$--dimensional objects, whereas our invariants measure {\em existence}
of $4$--dimensional things. In that sense, it might not come as a surprise
that there is an overlap between these theories. Note that the restriction
to trees is a well known feature if one wants {\em concordance invariants}
in the $3$--dimensional context, see \cite{CT2} or \cite{HM}.

\begin{figure}[ht!]
\centerline{\includegraphics[scale=.4]{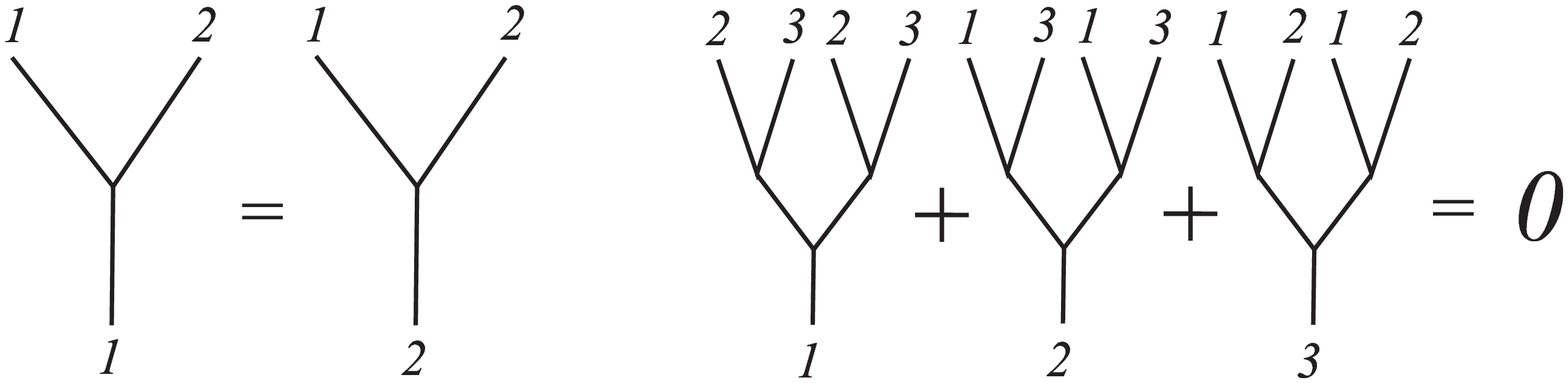}}
\caption{FR relations in order one and three (in a simply-connected
$4$--manifold).}
\label{FR-relations-fig}
\end{figure}

To make it possible that the intersection tree $\tau_n(\cW)$ only
depends on the $A_i$, it is in fact necessary to introduce two
more types of relations which correspond to changing the choices
of Whitney disks (for fixed choices of boundaries:
\begin{itemize}
\item The INT {\em interior or intersection} relations come from the
choice of the {\em interiors} of Whitney disks (which can be changed by
summing into any $2$--spheres). More generally, they measure
indeterminacies coming from certain lower order intersection trees for
Whitney towers on subsets of the $A_i$ together with other $2$--spheres. A
special case of these relations will be examined in detail in \cite{ST2}.
\item The FR {\em framing} relations are generated by certain $2$--torsion
elements which correspond to manipulations of the interiors of Whitney
disks that affect their normal framings. This will be described in
\cite{ST3} but see Figure~\ref{FR-relations-fig}.
\end{itemize}
The INT relations are more subtle in that they actually depend on
the ambient $4$--manifold $X$, rather than just on its fundamental
group. Both, INT and FR relations will not play a role in this
paper, however we will provide evidence supporting the following
conjecture by proving a closely related special case.
\begin{conj} \label{conj:well defined}
The intersection tree $\tau_n(\cW)\in \cT_n(\pi,m)/\mathrm{INT, FR}$
is independent of the choice of the Whitney tower $\cW$. In fact, it
only depends on the regular homotopy classes of the original maps $A_i$,
and should be written as $\tau_n(A_1,\dots,A_m)$.
\end{conj}
This result is well known in the Wall case, ie for $n=0$, and it was
proven in general for $n=1$ in \cite{ST} (and previously in the simply
connected case for $n=1$ in \cite{Ma} and \cite{FQ}).

The following result reflects the obstruction theoretic nature of the
intersection tree $\tau$.
\begin{thm} \label{thm:build-tower}
Let $A_i$ be properly immersed simply-connected surfaces in a
$4$--manifold, or connected surfaces in a simply-connected
$4$--manifold. If $\cW$ is an order--$n$ Whitney tower on the $A_i$ with
vanishing intersection tree $\tau_n(\cW)\in \cT_n(\pi,m)$, then there
is an order--$(n+1)$ Whitney tower on maps $A_i'$ which are regularly
homotopic (rel boundary) to $A_i$.
\end{thm}
Theorem~\ref{thm:build-tower} will be proved in
Section~\ref{sec:build-tower-proof}.

\subsection{Immersions with disjoint images} \label{sec:disjoint}
A special case of our invariant only counts those trees $t_p$
whose univalent labels are non-repeating, which means that the
number $m$ of spheres $A_i$ is two more than the order $n$ of the
intersection point $p$, $m=n+2$. Geometrically, one wants to
totally ignore {\em self}-intersections of the spheres $A_i$ and
in fact none of the (higher order analogues of) self-intersections
in the Whitney tower are paired up. This leads to the notion of a
{\em non-repeated Whitney tower} $\cW$ which has also a {\em
non-repeated intersection tree} $\lambda(\cW)$ that generalizes
the $\lambda$--invariant of Wall's intersection form. We shall
explain these notions in a different paper \cite{ST2} where we
also prove the following beautiful application of the theory.

\begin{thm} \label{thm:nonrepeating}
If the 2--spheres $A_1,\dots, A_{n+2}$ admit a non-repeated Whitney
tower $\cW$ of order~$n$, such that $\lambda(\cW)$ vanishes in
$\cT_n(\pi,n+2)$, then the homotopy classes (rel boundary) of the
$A_i$ can be represented by immersions with disjoint images.
\end{thm}

Again, this result was well known for $n=0$ (see eg \cite{Ko}),
and was proven for $n=1$ in \cite{ST} (and for trivial fundamental
group in \cite{Y}). In the special case discussed in the next
section, this result says that a link $L$ in $S^3$ has vanishing
non-repeating Milnor invariants if and only if it bounds disjoint
immersions of disks in $D^4$. In fact, this singular concordance
can then be improved to a {\em link homotopy} from $L$ to the
unlink (\cite{Go}, \cite{Gi}). This is Milnor's original theorem
\cite{M1}.

\subsection{Relation to Milnor invariants and the Kontsevich integral}
\label{sec:Milnor}
For a link $L \subset S^3$, there are unique homotopy classes (rel
boundary) $A_i\co D^2 \to D^4$ of immersions extending $L$. Therefore,
the previous discussion should apply to give link invariants via Whitney
towers. The {\em reduced} Kontsevich integral $Z^t(L)$ is the tree
part of the Kontsevich integral of $L$ and in \cite{HM} Habegger and
Masbaum have shown that the {\em first non-vanishing term} of $Z^t(L)-1$
carries exactly the same information as the first non-vanishing Milnor
invariants $\mu(L)$. These are the Milnor invariants with repeating
indices, also denoted $\bar\mu$--invariants \cite{M2}. We shall not
make this distinction and we consider only the ``first non-vanishing''
invariants. In the general case one needs to consider {\em string}
links \cite{HM}.

Denote by $K_n(L)$ the order--$n$ term of $Z^t(L)-1$. Now observe
that $K_n(L)$ takes values exactly in $\cT_n(m)\otimes\Q$, where
$m$ is the number of components of $L$ and the {\em order}~$n$ is
the number of trivalent vertices. Here the relations in $\cT(m)$
simplify dramatically because $\pi_1(D^4)=0=\pi_2(D^4)$ and in
fact they reduce to exactly the AS and IHX relations used in the
usual definition of the Kontsevich integral. We note that the most
commonly used degree in papers on the Kontsevich integral is one
half the total number of vertices. For unitrivalent trees, this
degree is one more than the number of trivalent vertices, ie one
{\em more} than the order that we are using here.

For an oriented link $L\subset S^3$, consider the following four statements.
\begin{enumerate}
\item $L$ bounds a Whitney tower of order~$n$ in $D^4$.
\item $L$ bounds disjointly embedded framed gropes of class~$(n+1)$ in $D^4$.
\item $L$ has vanishing $\mu$--invariants of length~$\leq (n+1)$.
\item All terms in $Z^t(L)-1$ having order~$\leq (n-1)$ vanish.
\end{enumerate}
Then $\mathrm{(i)}$ is equivalent to $\mathrm{(ii)}$ by \cite{S},
$\mathrm{(iii)}$ is equivalent to $\mathrm{(iv)}$ by \cite{HM}, and
$\mathrm{(ii)}$ implies $\mathrm{(iii)}$ by \cite{KT}.

The following theorem gives the relation between the Kontsevich integral
and our intersection tree $\tau$ in the context of the above results.
\begin{thm} \label{thm:Milnor}
If $L$ bounds a Whitney tower $\cW$ of order~$n$ in $D^4$, then
$$K_n(L)=\tau_n(\cW) \quad \in \cT_n(m)\otimes\Q$$
which shows that rationally, $\tau_n(L):=\tau_n(\cW)$ only depends on
(the concordance class of) $L$ and can be used to calculate the first
non-vanishing terms of the reduced Kontsevich integral as well as the
Milnor invariants.
\end{thm}

\begin{rem} In \cite{ST3} we shall explain a direct geometric relation between our intersection trees
and Milnor's invariants, completely avoiding the Kontsevich integral.
\end{rem}

\begin{rem} \label{rem:torsionfree}
In the nonrepeating case, the groups $\cT_n(n+2)$ are torsionfree,
and hence tensoring with $\Q$ does not lose any information. This
implies our above Conjecture~\ref{conj:well defined} for this very
special case (since the FR and INT relations are trivial). By
results in \cite{Ko}, Theorem~\ref{thm:Milnor} also implies the
conjecture for the $2$--spheres in the simply connected
$4$--manifold formed by attaching 0--framed $2$--handles to the
$4$--ball along $L$ in the nonrepeating case (or rationally in the
repeating case). It is not unreasonable to believe that the groups
$\cT_{2n}(m)$ are also torsionfree (with repeated labels allowed).
Note that $\cT_1(1) \cong \Z/2$ which corresponds exactly to the
Arf invariant of a knot (see \cite{Ma}, \cite{S2}, \cite{ST}) and
hence shows that statement (iv) does {\em not} imply (i) in the
above theorem. In general, the FR relations are non-trivial for
odd orders as will be explained in \cite{ST3}; see
Figure~\ref{FR-relations-fig} for an example.
\end{rem}

\section{Whitney towers and intersection trees}\label{W-tower-tau-sec}
The goal of this section is to define the $n$th-order intersection
tree $\tau_n(\cW)$ of an order--$n$ Whitney tower $\cW$ in an oriented
$4$--manifold $X$. After giving the precise definition of a Whitney
tower $\cW$, an indexing of the surfaces in $\cW$ is given in terms
of bracketings and rooted trees which are labelled, oriented and then
decorated by elements of the fundamental group $\pi :=\pi_1X$. The
unrooted decorated tree $t_p$ associated to an intersection point $p$
in $\cW$ then corresponds to a pairing of the rooted trees associated to
the intersecting surfaces. Finally, $\tau_n(\cW)$ is defined as a signed
sum of the $t_p$ in the group $\cT_n(\pi,m)$, see Section~\ref{subsec:T}.

\subsection{Whitney towers}\label{W-tower-sub-sec}
We assume our $4$--manifolds are oriented and equipped with a
basepoint. The reader is referred to \cite{FQ} for details on immersed
surfaces in $4$--manifolds, including {\em Whitney moves} and (Casson)
{\em finger moves}. For more on Whitney towers see \cite{CST}, \cite{S},
\cite{S2}.

\begin{defn}\label{w-tower-defn}\mbox{}
\begin{itemize}
\item A {\em surface of order 0} in a 4--manifold $X$
is a properly immersed surface (boundary embedded in the boundary
of $X$ and interior immersed in the interior of $X$). A {\em
Whitney tower of order 0} in $X$ is a collection of order--0
surfaces.
\item The {\em order of a (transverse) intersection point} between a
surface of order $n_1$ and a surface of order $n_2$ is $n_1+n_2$.
\item The {\em order of a Whitney disk} is $n+1$ if it pairs intersection
points of order $n$.
\item For $n\geq 0$,
a {\em Whitney tower $\cW$ of order $n+1$}  is a Whitney tower of
order $n$ together with Whitney disks pairing all order--$n$
intersection points of $\cW$. These top order disks are allowed
to intersect each other as well as lower order surfaces.
\end{itemize}

The Whitney disks in a Whitney tower are required to be {\em
framed} (\cite{FQ}) and have disjointly embedded boundaries.
Intersections in surface interiors are assumed to be transverse. A
Whitney tower is {\em oriented} if all its surfaces (order--$0$
surfaces {\em and} Whitney disks) are oriented. A {\em based}
Whitney tower includes a chosen basepoint on each surface
(including Whitney disks) together with a {\em whisker} (arc) for
each surface connecting the chosen basepoints to the basepoint of
the ambient $4$--manifold.
\end{defn}

Some further terminology: If $\cW$ is an order--$n$ Whitney tower
containing $A_i$ as its order--$0$ surfaces then the $A_i$ are said
to {\em admit} an order--$n$ Whitney tower and we say that $\cW$ is
a Whitney tower {\em on} the $A_i$.

\subsection{Rooted trees and brackets}\label{trees-brackets}
Non-associative ordered {\em bracketings} of elements from some
index set correspond to {\em rooted labelled vertex-oriented
unitrivalent trees} as follows. Here {\em rooted} means ``having a
preferred univalent vertex'' (the {\em root}), {\em labelled}
means that each non-root univalent vertex is labelled by an
element from the index set and {\em vertex-oriented} means that
each trivalent vertex is equipped with a cyclic ordering of its
incident edges. The {\em order} of a tree is the number of
trivalent vertices.

\begin{figure}[ht!]
\centerline{\includegraphics[scale=.5]{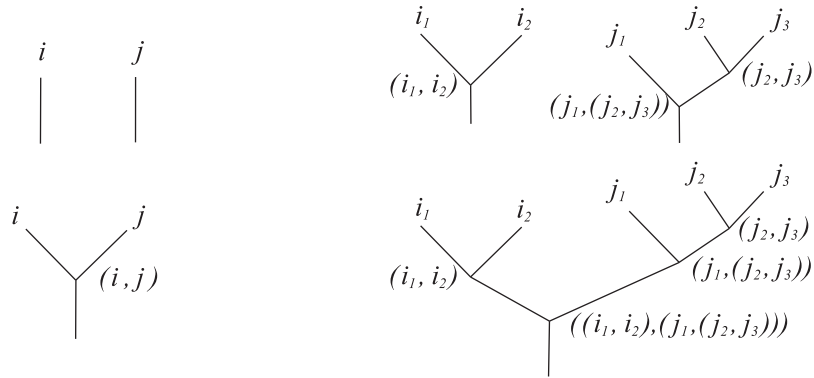}}
\caption{Rooted trees $t(i)$ and $t(j)$ (upper left) and their rooted product
$t(i,j)=t(i)\ast t(j)$ (lower left);
$t(i_1,i_2)$ and $t(j_1,(j_2,j_3))$ (upper right) and their rooted product
$t((i_1,i_2),(j_1,(j_2,j_3)))=t(i_1,i_2)\ast t(j_1,(j_2,j_3))$
(lower right). In this figure all trivalent orientations are clockwise
in the plane.}
\label{bracket-treesA-fig}
\end{figure}

A bracketing $(i)$ of a singleton element $i$ from the index set
corresponds to the rooted order--0 tree $t(i)$ consisting of a single
edge with one vertex labelled by $i$ and the other vertex designated as
the root. A bracketing $(I,J)$ of brackets $I$ and $J$ corresponds to
the {\em rooted product} $t(I,J):=t(I) \ast t(J)$ of the trees $t(I)$
and $t(J)$ which identifies together the roots of $t(I)$ and $t(J)$
to a single vertex and ``sprouts'' a new rooted edge at this vertex
(Figure~\ref{bracket-treesA-fig}) with the cyclic order at the new
trivalent vertex given by taking the edges coming from $I$, $J$ and the
root in that order.

Thus, the non-root univalent vertices of the tree $t(I)$
associated to a bracket $I$ are labelled by elements from the
index set and the trivalent vertices correspond to sub-bracketings
of $I$, with the trivalent vertex adjacent to the root
corresponding to $I$.

\begin{rem}\label{rem:w-tree-f-move}
The rooted product $\ast$ can be ``realized'' geometrically by a
finger-move: Pushing a Whitney disk $W_I$ through another Whitney
disk $W_J$ creates $W_{(I,J)}$ with $t(W_{(I,J)})=t(W_I)\ast
t(W_J)$.
\end{rem}

This remark uses the upcoming assignment of a rooted tree $t(W)$
to a Whitney disk $W$ inside a Whitney tower $\cW$. In the easiest
version, one starts with a root for $W$ and then introduces one
branching (trivalent vertex) while reading off which two sheets of
$\cW$ are paired by $W$. Then one continues with the same
procedure for the two sheets to inductively obtain $t(W)$. In the
next section we shall make this procedure precise, and in fact
explain directly how orientations on the Whitney disks lead to
vertex-orientations of the corresponding trees.

\subsection{Rooted trees for oriented Whitney towers}
Let $\cW$ be an oriented Whitney tower on order--$0$ surfaces $A_i$
for $i=1,2,\ldots,m$. The orientations on the surfaces in $\cW$
set up an indexing of the surfaces in $\cW$ by bracketings $I$
from $\{1,2,\ldots,m\}$ and their corresponding rooted vertex
oriented unitrivalent $m$--labelled trees $t(I)$
(\ref{trees-brackets}) via the following conventions:

A bracketing $(i)$ of a singleton element $i$ from the index set
and the corresponding rooted order--0 tree $t(A_i):=t(i)$ are
associated to each order--$0$ surface $A_i$. The bracket $(I,J)$
and the corresponding tree $t(W_{(I,J)}):=t(I,J)$ are associated
to a Whitney disk $W_{(I,J)}$, pairing intersections between $W_I$
and $W_J$, with the ordering of the components $I$ and $J$ in the
associated bracket $(I,J)$ chosen so that the orientation of
$W_{(I,J)}$ is the same as that given by orienting its boundary
$\partial W_{(I,J)}$ from the negative intersection point to the
positive intersection point {\em first} along $W_I$ {\em then}
back along $W_J$ to the negative intersection point, together with
a second inward pointing tangent vector.

We use brackets as subscripts to index surfaces in $\cW$, writing
$A_i$ for an order--$0$ surface (dropping the brackets around the
singleton $i$) and $W_{(i,j)}$ for a first-order Whitney disk that
pairs intersections between $A_i$ and $A_j$, etc.. When writing
$W_{(I,J)}$ for a Whitney disk pairing intersections between $W_I$
and $W_J$, the understanding is that if a bracket $I$ is just a
singleton $(i)$ then the surface $W_I=W_{(i)}$ is just the
order--$0$ surface $A_i$. In general, the order of $W_I$ is equal
to the order of (ie the number of trivalent vertices of)
$t(W_I)$.

\begin{figure}[ht!]
\centerline{\includegraphics[scale=.5]{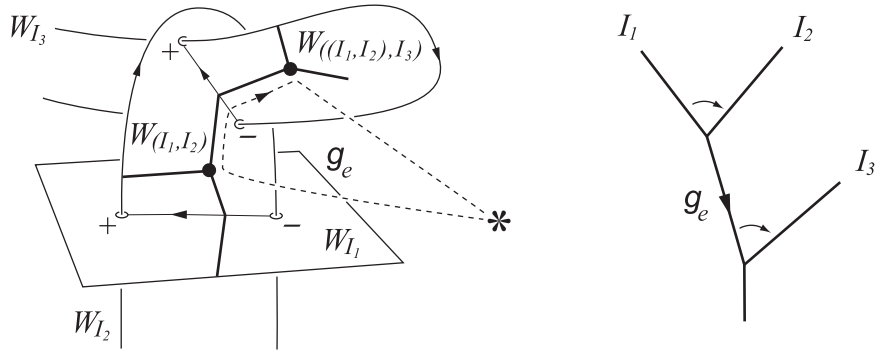}}
\caption{A Whitney disk $W_{((I_1,I_2),I_3)}$ and its
associated tree $t(W_{((I_1,I_2),I_3)})$ shown (left) as a
subset of the Whitney tower and (right) as an abstract
rooted tree. The boundaries of the Whitney disks are
oriented according to our bracket-orientation conventions
using the indicated signs of the intersection points. The
dashed path indicates a sheet-changing loop (based at the
basepoint of the ambient $4$--manifold $X$) which
determines the element $g_e\in\pi_1X$ decorating the
corresponding oriented edge as described in
\ref{decorating-tree-subsec}.} \label{W-disk-treeC-OR-fig}
\end{figure}

It will be helpful to consider each tree $t(W_I)$ as a subset of
$\cW$: Assuming that $\cW$ is based (Definition~\ref{w-tower-defn}),
map the vertices (other than the root) of $t(W_I)$ to the basepoints
of the surfaces whose indices are contained as sub-brackets of $I$
and map the edges (other than the edge adjacent to the root) of
$t(W_I)$ to sheet-changing paths between basepoints, as illustrated in
Figure~\ref{W-disk-treeC-OR-fig} (disregarding, for the moment, the dotted
loop which will be explained in \ref{decorating-tree-subsec}). Then embed
the root and its edge anywhere in the {\em negative corner} of $W_I$
(see next paragraph).

It can be arranged that this mapping of $t(W_I)$ into $\cW$ has the
property that the trivalent orientations of $t(W_I)$ are induced by
the orientations of the corresponding Whitney disks: Note that the pair
of edges which pass from a trivalent vertex down into the lower order
surfaces paired by a Whitney disk determine a ``corner'' of the Whitney
disk which does not contain the other edge of the trivalent vertex. If
this corner contains the {\em positive} intersection point paired by
the Whitney disk, then the vertex orientation and the Whitney disk
orientation agree. Our figures are drawn to satisfy this convention.

\subsection{Orientation choices on Whitney disks}
\label{subsec:orientations}
Via our bracket-orientation convention, changing the orientation on a
Whitney disk $W_{(I,J)}$ changes its tree from $t(W_{(I,J)})=t(I,J)$ to
$t(W_{(J,I)})=t(J,I)$, ie changes the cyclic orientation of the associated
trivalent vertex. In addition, changing the orientation of a {\em single}
lower order Whitney disk $W_K$ corresponding to a trivalent vertex of
$t(W_{(I,J)})$ (so $K$ is a sub-bracket of $(I,J)$, with $K\neq (I,J)$)
changes the cyclic orientations at exactly {\em two} trivalent vertices of
$t(W_{(I,J)})$: the one corresponding to $W_K$ and the adjacent one which
corresponds to a Whitney disk pairing intersections between $W_K$ and some
other surface. This is because changing the orientation of $W_K$ reverses
the signs of the intersection points between $W_K$ and anything else.

\subsection{Decorated trees for Whitney towers}
\label{decorating-tree-subsec}
Let $t(W_I)$ be the (oriented labelled rooted) tree associated to
a Whitney disk $W_I$ in an oriented based Whitney tower $\cW$ in
a $4$--manifold $X$. Thinking of $t(W_I)$ as a subset of $\cW$ as
described above, any edge $e$ of $t(W_I)$, other than the root-edge,
corresponds to a sheet-changing path connecting the basepoints
of adjacent surfaces in $\cW$. For a chosen orientation of $e$,
this path together with the whiskers on  the adjacent surfaces form
an oriented loop which determines an element $g_e$ of $\pi:=\pi_1X$
(Figure~\ref{W-disk-treeC-OR-fig}). Fixing (arbitrarily) orientations for
all the (non-root) edges in $t(W_I)$ and labelling each oriented edge
with an element of $\pi$ in this way yields the {\em decorated rooted
tree} associated to $W_I$ (which will still be denoted by $t(W_I)$). Note
that switching the orientation of $e$ changes $g_e$ to $g_e^{-1}$ which
explains the OR orientation reversal relation mentioned in \ref{tau-intro}
and shown in Figure~\ref{Relations-fig}. (Since we are working in an
orientable $4$--manifold, $\omega_1(g_e)$ is trivial.) Also, changing the
choice of whisker on a Whitney disk has the effect of left multiplication
on the group elements associated to the three edges adjacent to and
oriented away from the trivalent vertex corresponding to the Whitney
disk accounting for the HOL relation.

When decorations are understood, we will also denote a decorated tree
by $t(I)$ where the underlying tree corresponds to the bracket $I$.

\subsection{Decorated trees for intersection points}
\label{int-point-trees}
If $p$ is a transverse intersection point between $W_I$ and $W_J$
in $\cW$ then the {\em decorated tree} $t_p$ associated to $p$ is
defined as follows. Identify the roots of the decorated trees $t(W_I)$
and $t(W_J)$ to a single (non-vertex) point. The two edges that were
adjacent to the roots of $t(W_I)$ and $t(W_J)$ now form a single edge
$e_p$. Chose an orientation of $e_p$ and decorate $e_p$ by the element
of $\pi$ determined by the whiskers on $W_I$ and $W_J$ together with a
path connecting the basepoints of $W_I$ and $W_J$ that changes sheets
only at $p$ with the orientation induced by $e_p$.

\begin{figure}[ht!]
\centerline{\includegraphics[scale=.5]{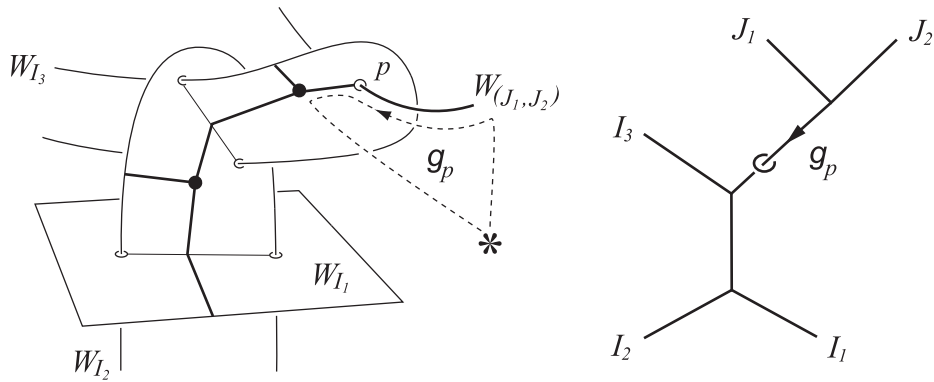}}
\caption{The punctured tree $t^{\circ}_p$ associated to an
intersection point $p\in W_I\cap W_J$ (for
$I=((I_1,I_2),I_3)$ and $J=(J_1,J_2)$) shown as a subset
of the Whitney tower and as an abstract labelled
(punctured) tree. Decorations other than $g_p$ are
suppressed and the sheet-changing loop that determines
$g_p$ is indicated by the dashed path.}
\label{W-disk-int-point-treeD-OR-fig}
\end{figure}

Thus, the decorated tree $t_p$ is unrooted and every edge of $t_p$
is oriented and decorated with an element of $\pi$. Note that the
order of $p$ is equal to the order of $t_p$ (the number of
trivalent vertices).

The mappings of $t(W_I)$ and $t(W_J)$ into $\cW$ give rise to a
mapping of $t_p$ into $\cW$: Just map the root vertices of $W_I$
and $W_J$ to $p$ and the adjacent edges become a sheet-changing
path between the basepoints of $W_I$ and $W_J$
(Figure~\ref{W-disk-int-point-treeD-OR-fig}). This mapping is an
embedding of $t_p$ into $\cW$ if all the Whitney disks ``beneath''
$W_I$ and $W_J$ (corresponding to sub-brackets of $I$ and $J$) are
distinct.

We will sometimes keep track of the edge of $t_p$ that corresponds
to $p$ by marking that edge with a small linking circle as in
Figure~\ref{W-disk-int-point-treeD-OR-fig}; such a {\em punctured tree}
will be denoted by $t^{\circ}_p$.

\begin{figure}[ht!]
\centerline{\includegraphics[scale=.5]{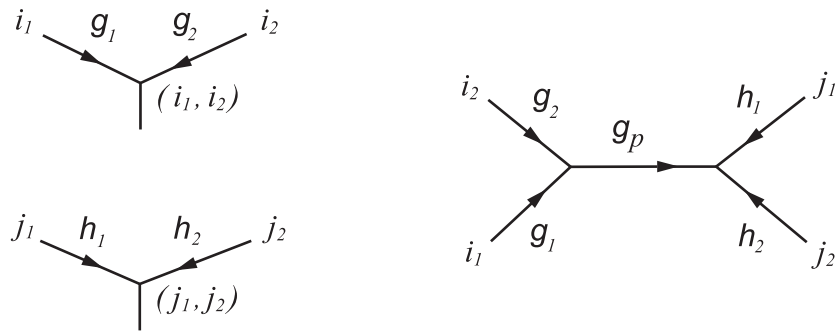}}
\caption{A pair of decorated rooted trees $t(I)$ and $t(J)$
corresponding to order--1 Whitney disks $W_{I}$ and
$W_{J}$ with $I=(i_1,i_2)$ and $J=(j_1,j_2)$ (left), and the inner product
$t_p=t(W_I)\cdot_{g_p} t(W_J)=t(I)\cdot_{g_p} t(J)$ associated to
an order--2 intersection point $p\in W_{I}\cap W_{J}$ (right).}
\label{inner-product-treesC-fig}
\end{figure}

It will be convenient to formalize the above description of the
(unrooted) decorated tree $t_p$ as a pairing (over the group
$\pi$) of rooted decorated trees: Given a pair $t(I)$ and $t(J)$
of rooted decorated trees and an element $g\in\pi$, define the
{\em inner product} $t(I)\cdot_g t(J)$ to be the unrooted
decorated tree gotten by identifying together the root vertices of
$t(I)$ and $t(J)$ to a single (non-vertex) point in an edge
labelled by $g$ as illustrated in
Figure~\ref{inner-product-treesC-fig}. Thus, in this notation we
have $t_p:=t(W_I)\cdot_{g_p} t(W_J)$ for $p\in W_I\cap W_J$ as
just described above.

\subsection{The antisymmetry AS relation}\label{subsec:AS}
If a Whitney tower $\cW$ is oriented then there is one more piece
of information that we need to keep track of: the sign $\e_p$ of
an unpaired intersection point
$$p\in W_I\cap W_J\subset \cW.$$
$\e_p$ is computed, in the usual way, by comparing the orientation
determined by $W_I$ and $W_J$ at $p$ with the orientation of the
ambient $4$--manifold $X$ at $p$.

Changing the orientation on the Whitney disk $W_I$ changes the
{\em signed tree} $\epsilon_p\cdot t_p$ by the AS antisymmetry
relation mentioned in \ref{tau-intro}: The cyclic orientation of
the vertex corresponding to $W_I$ in $t_p$ is switched and so is
the sign $\epsilon_p$ of the intersection with $W_J$. Moreover,
changing the orientation of a single Whitney disk, other than
$W_I$ or $W_J$, preserves the sign $\epsilon_p$ and changes the
cyclic orientations at {\em two} trivalent vertices of $t_p$, as
pointed out above in Section~\ref{subsec:orientations}.
Consequently, working modulo the AS relation makes the signed tree
$\epsilon_p\cdot t_p$ independent of the choices of orientations
for the Whitney disks in $\cW$.

The dependence on orientations for the original sheets $A_i$
remains: changing the orientation of one $A_i$ introduces an
additional sign into $\epsilon_p\cdot t_p$ if $t_p$ has an odd
number of $i$--labelled vertices.

\subsection{The intersection tree $\tau_n(\cW)$}\label{subsec:T}
\begin{figure}[ht!]
\centerline{\includegraphics[scale=0.9]{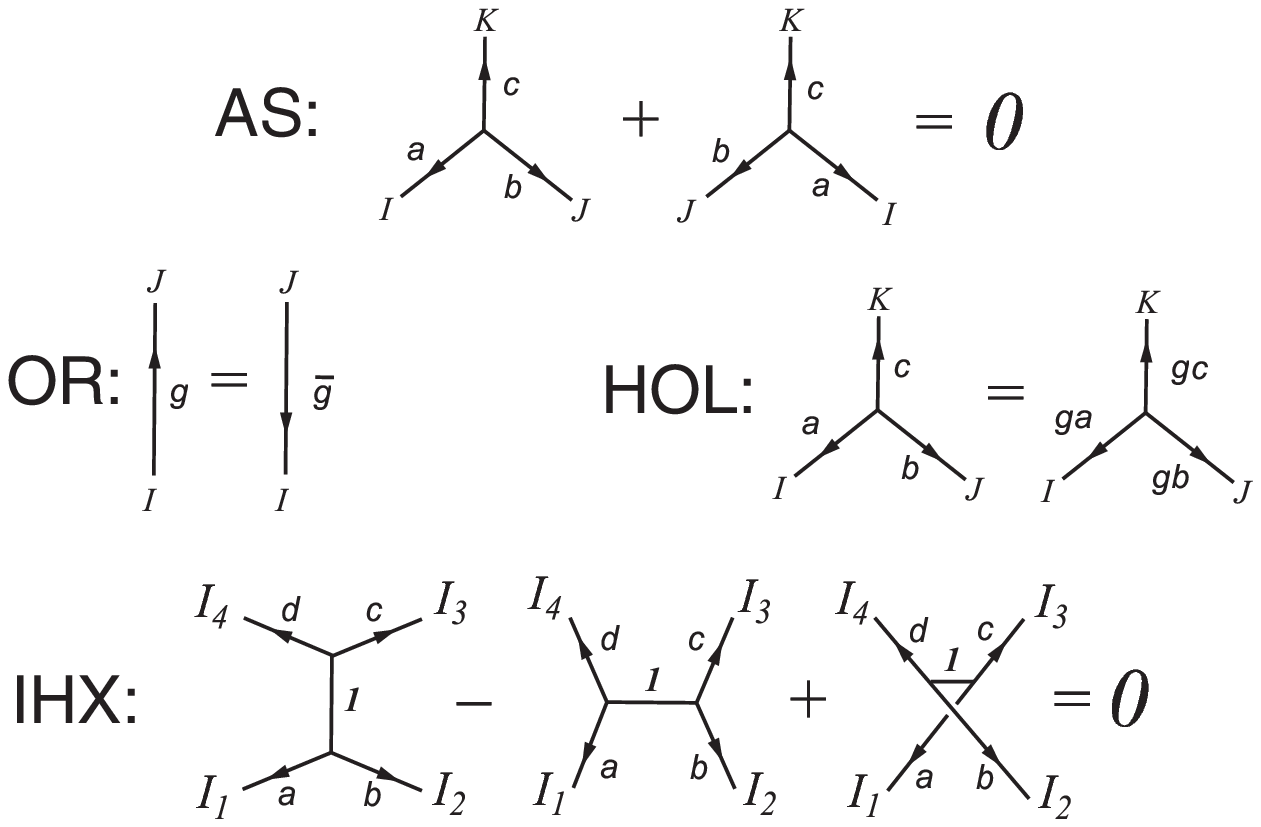}}
\caption{The AS, OR, HOL and IHX relations in $\cT_n(\pi,m)$ for
$a$, $b$, $c$, $d$, $1$ and $g$ in $\pi$ with $\overline{g}=g^{-1}$. All
trivalent orientations are induced from a fixed orientation of the plane.}
\label{Relations-fig}
\end{figure}

We would next like to add up the unpaired intersection points of a
given Whitney tower in some algebraic structure. For that purpose,
let $\cT_n(\pi,m)$ denote the abelian group generated by
(isomorphism classes of) decorated trees of order~$n$ modulo the
relations shown in Figure~\ref{Relations-fig}. That is, each
generator is an (unrooted) unitrivalent tree having
\begin{itemize}
\item  $n$ cyclically oriented trivalent vertices,
\item $n+2$ univalent vertices labelled by elements of $\{1,\dots,m\}$, and
\item $2n+1$ oriented edges decorated by elements of $\pi$.
\end{itemize}

\begin{defn}\label{int-tree-defn}
Let $\cW$ be an order--$n$ Whitney tower on properly immersed
simply-connected oriented surfaces $A_1,\dots,A_m$ in a
$4$--manifold $X$. (In fact, the $A_i$ only need to be
$\pi_1$--null, see \cite{FQ}.) Define the  $n$th-order {\em
intersection tree} of $\cW$ by
$$\tau_n(\cW):= \sum_p \e_p\cdot t_p \quad \in \quad \cT_n(\pi,m)$$
where the sum is over all order--$n$ intersection points $p$ in $\cW$.
\end{defn}

As explained above, the AS relations make sure that $\tau_n(\cW)$
actually does not depend on the choice of orientations for the
Whitney disks. Similarly, the HOL and OR relations make sure that
$\tau_n(\cW)$ does not depend on the choice of whiskers, or edge
orientations. In other words, $\tau_n(\cW)$ is defined by first
choosing whiskers and orientations (on edges and Whitney disks)
and then proving independence of these choices.

\begin{rem} Using the HOL relation or, more concretely, by choosing
the whiskers on the Whitney disks appropriately, one can normalize
the trees $t_p$ so that all interior edges and one univalent edge
are decorated with the trivial group element $1\in\pi$. Thus, one
can interpret $\tau_n(\cW)$ as living in a quotient of the
integral group ring of the $(n+1)$--fold product of $\pi$.

By slightly refining our notation, signs can be associated
formally to all tree edges and the edge decorations can be
extended linearly to elements of the group ring $\Z [\pi ]$
(compare \cite{GK}, \cite{GL}). Similarly, one can extend the
labels on the univalent vertices to the free abelian group on
$\{1,\dots,m\}$.
\end{rem}

\section{Proof of Theorem~\ref{thm:build-tower}}\label{sec:build-tower-proof}
Our proof of Theorem~\ref{thm:build-tower} will be constructive in
the sense that we describe how to build the next order Whitney
tower by geometrically realizing all the relations in
$\cT_n(\pi,m)$. However, it should be mentioned that since the
groups $\cT_n(\pi,m)$ do not in general have a canonical basis we
are sidestepping the ``word problem'' in $\cT_n(\pi,m)$.  The main
construction (Lemma~\ref{transfer-lemma}) of the proof shows how
to exchange {\em algebraic cancellation} of pairs of intersection
points for {\em geometric cancellation} (by Whitney disks) in the
case that the intersection points are {\em simple} (have certain
standard right- or left-normed trees,
\ref{simple-int-transfer-lem-subsec}). This algebraic cancellation
occurs in the lift $\widehat{\cT}$ of $\cT$ which forgets the IHX
relation. The general case is then reduced to this case using
geometric IHX constructions from \cite{CST} and \cite{S} to show
that an order--$n$ Whitney tower $\cW$ with
$\tau_n(\cW)=0\in\cT_n(\pi,m)$ can be modified so that all
order--$n$ intersections come in simple algebraically-cancelling
pairs.

To simplify the exposition and highlight the combinatorial structure
of Whitney towers, we will emphasize the simply-connected case, often
dropping the group $\pi$ from notation. Refining the constructions to
cover the general case for the most part only requires checking that
whiskers can be (re)-chosen appropriately. At a first reading it doesn't
hurt to ignore group elements entirely and only the simply-connected
version of Theorem~\ref{thm:build-tower} will be used later in the proof
of Theorem~\ref{thm:Milnor}.

We begin with some notation and lemmas. All Whitney towers are assumed
oriented, labelled and based.

\subsection{Geometric intersection trees for Whitney towers}
\label{geo-int-tree-w-tower}
For an (oriented, labelled, based) Whitney tower $\cW$ define $t_n(\cW)$,
the ($n$th-order, oriented) {\em geometric intersection tree} of $\cW$,
to be the disjoint union of signed (decorated) trees
$$t_n(\cW):=\amalg_p \e_p\cdot t_p$$
over all unpaired order--$n$ intersection points $p \in \cW$. (An unsigned
version of $t_n(\cW)$ was defined for unoriented Whitney towers in
\cite{S}.) The next two pairs of definitions and lemmas will illustrate
how $t_n(\cW)$ captures the essential geometric structure of $\cW$.

\subsection{Split subtowers}\label{split-subtowers}
The Whitney disks in an arbitrary Whitney tower may have multiple self-intersections and
intersections with other surfaces. However, it is not difficult to modify an arbitrary Whitney
tower so that each Whitney disk is embedded and contains either a single Whitney arc or unpaired
intersection point (Lemma~\ref{split-tower-lem} below). This is best expressed using the notion of
{\em split subtowers} and splitting a Whitney tower into split subtowers will serve to simplify
geometric constructions and combinatorial arguments.

The purpose of constructing a Whitney tower is to provide
information on the homotopy classes (rel boundary) of its
order--$0$ surfaces. However, when describing and manipulating {\em
subsets} of a Whitney tower it is natural to consider {\em
sub}towers on sheets of surfaces which are not {\em properly}
immersed:
\begin{defn}\label{subtower-defi}
A {\em subtower} is a Whitney tower except that the boundaries of
the immersed order--$0$ surfaces in a subtower are allowed to lie
in the interior of the $4$--manifold (instead of being required to
lie in the boundary). The boundaries of the order--$0$ surfaces in
a subtower are still required to be embedded. The notions of {\em
order} for intersection points and Whitney disks are the same as
in Definition~\ref{w-tower-defn}.
\end{defn}
In this paper we will only be concerned with subtowers whose order--0
surfaces are sheets in the order--0 surfaces of an actual Whitney tower. In
this case, the surfaces of the subtower inherit the same orientations
and indexing by brackets as the Whitney tower. Thus, the association of
decorated trees to surfaces and intersection points is also the same.

\begin{defn}\label{split-subtower-defi}
A subtower $\cW_p$ is {\em split} if it satisfies all of the following:
\begin{enumerate}
\item $\cW_p$ contains a single unpaired intersection point $p$,
\item the order--$0$ surfaces of $\cW_p$ are all embedded $2$--disks,
\item the Whitney disks of $\cW_p$ are all embedded,
\item the interior of any surface in $\cW_p$ either contains $p$ or
contains a single Whitney arc of a Whitney disk in $\cW_p$,
\item $\cW_p$ is connected (as a $2$--complex in the $4$--manifold).
\end{enumerate}
Moreover, a Whitney tower $\cW$ is called {\em split} if all the
unpaired intersection points of $\cW$ are contained in disjoint
split subtowers on sheets of the order--$0$ surfaces of $\cW$.
\end{defn}

Note that a normal thickening of a split subtower $\cW_p$ in the
ambient $4$--manifold is just the $4$--disk $D^4$ which is  a
regular neighborhood of the embedded tree $t_p$ associated to the
unpaired intersection point $p$.

\begin{figure}[ht!]
\centerline{\includegraphics[scale=.5]{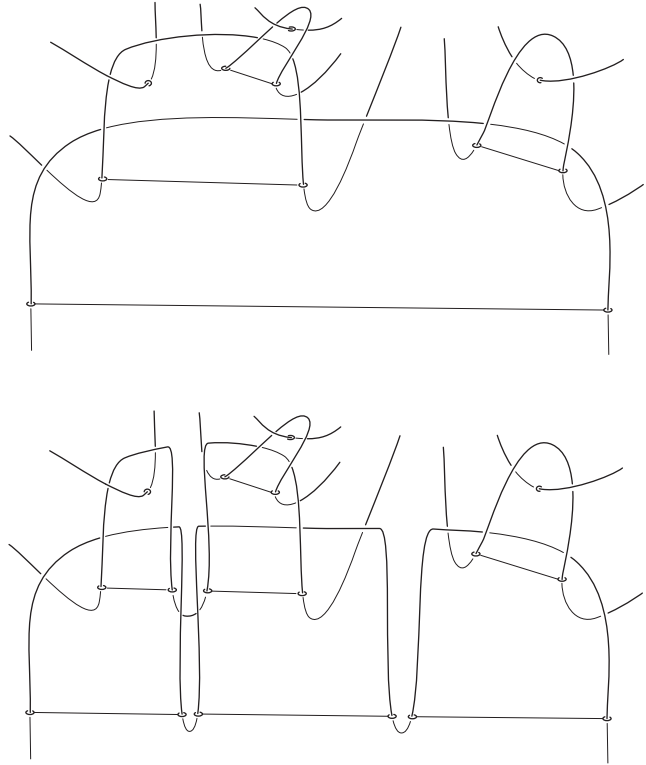}}
\caption{Part of a Whitney tower before (top) and after (bottom) applying
the splitting procedure described in the proof of Lemma~\ref{split-tower-lem}.}
\label{splitting-and-split-tower-fig}
\end{figure}

\subsection{Split Whitney towers}\label{split-w-tower-subsec}
The splitting of a Whitney tower into split subtowers described in the
following lemma is analogous to Krushkal's splitting of a grope into
genus one gropes \cite{Kr}.
\begin{lem}\label{split-tower-lem}
Let $\cW$ be a Whitney tower on order--0 surfaces $A_i$. Then there exists
a split Whitney tower $\cW_{\mathrm{split}}$ contained in any regular
neighborhood of $\cW$ such that:
\begin{enumerate}
\item The order--0 surfaces $A'_i$ of $\cW_{\mathrm{split}}$ only differ
from the $A_i$ by finger moves.
\item The geometric intersection trees $t(\cW)$ and
$t(\cW_{\mathrm{split}})$ are isomorphic.
\end {enumerate}
\end{lem}
The isomorphism in item~(ii) includes decorations and signs.

\begin{proof}
Starting with the highest-order Whitney disks of $\cW$, apply finger moves
as indicated in Figure~\ref{splitting-and-split-tower-fig}. Working
down through the lower-order Whitney disks yields the desired
$\cW_{\text{split}}$. Choosing whiskers and orientations appropriately
for the new Whitney disks preserves the decorations on the trees
associated to the unpaired intersection points.
\end{proof}

An advantage of splitting a Whitney tower is that the geometric
intersection tree sits as an {\em embedded} subset (\ref{int-point-trees})
and all the singularities of the split Whitney tower are contained in
disjointly embedded $4$--balls, each of which is a regular neighborhood
of an intersection point tree. In this sense the decomposition of a
Whitney tower into split subtowers corresponds to the idea that the trees
associated to the unpaired intersection points capture the essential
structure of a Whitney tower. The next lemma can be interpreted as
justifying that this essential structure is indeed captured by the {\em
un}-punctured trees rather than the punctured trees in the sense that an
unpaired intersection point (corresponding to a punctured edge) can be
``moved'' to any other edge of its tree.

\begin{figure}[ht!]
\centerline{\includegraphics[scale=.6]{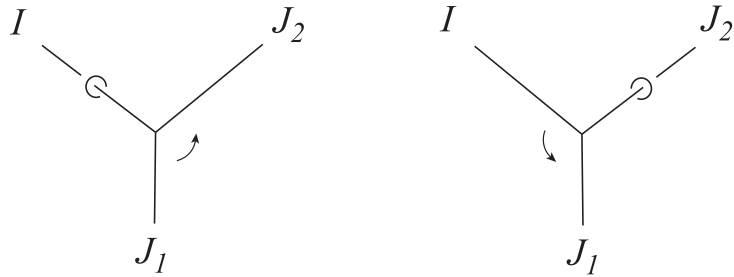}}
\caption{A local picture of the tree associated to the split
subtower $\cW$
before (left) and $\cW'$ after (right) the Whitney move in the proof
of Lemma~\ref{subtower-lemma} illustrated in
Figure~\ref{tree-move1-OR-fig} and Figure~\ref{tree-move2-OR-fig}.}
\label{tree-move-trees-fig}
\end{figure}

\begin{lem}\label{subtower-lemma}
Let $\cW\subset X$ be a split subtower on order--$0$ sheets $s_i$
with unpaired intersection point $p=W_I\cap W_J\subset \cW$.
Denote by $\nu(\cW)$ a normal thickening of $\cW$ in $X$ so that
$\partial s_i\subset\partial\nu(\cW)\subset\nu(\cW)\cong D^4$. If
$I'$ and $J'$ are any brackets such that the decorated trees
$t(I')\cdot t(J')=t_p=t(I)\cdot t(J)$, then after a homotopy (rel
$\partial$) of the $s_i$ in $\nu(\cW)$ the $s_i$ admit a split
subtower $\cW'\subset\nu(\cW)$ with single unpaired intersection
point $p'=W_{I'}\cap W_{J'}\subset \cW'$ such that
$\epsilon_{p'}\cdot t_{p'}=\epsilon_p\cdot t_p$.
\end{lem}

\begin{figure}[ht!]
\centerline{\includegraphics[scale=.55]{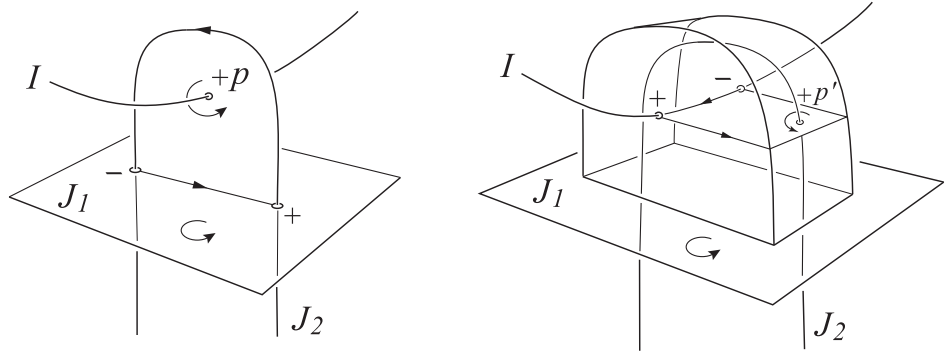}}
\caption{The unpaired intersection point $p=W_I\cap W_J$ in the split
subtower $\cW$ of Lemma~\ref{subtower-lemma} (left), and the unpaired
intersection point $p'=W_{I'}\cap W_{J'}$ in $\cW'$ after the Whitney move (right).
Signs and orientations are indicated for the case $\e_p=+$, with brackets
corresponding to the trivalent orientations in Figure~\ref{tree-move-trees-fig}.}
\label{tree-move1-OR-fig}
\end{figure}
\begin{figure}[ht!]
\centerline{\includegraphics[scale=.55]{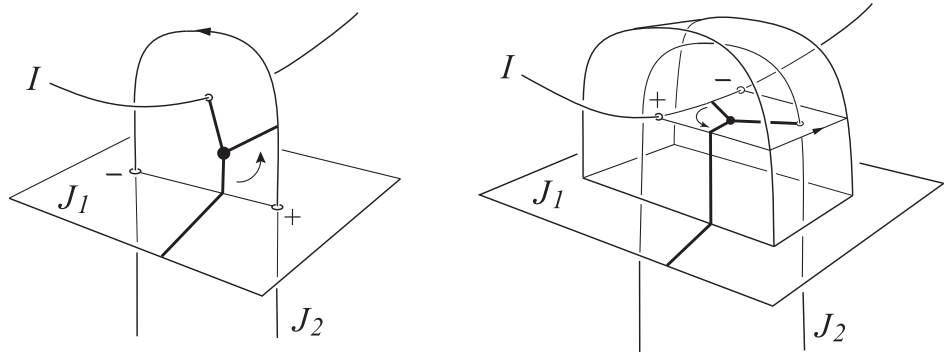}}
\caption{This figure shows that the oriented punctured trees associated
to $p$ and $p'$ in Figure~\ref{tree-move1-OR-fig} differ as indicated in
Figure~\ref{tree-move-trees-fig}.}
\label{tree-move2-OR-fig}
\end{figure}
\begin{proof}(of Lemma~\ref{subtower-lemma})
It is enough to show that the puncture in $t^{\circ}_p$ can be
``moved'' to either {\em adjacent} edge, since by iterating it can
be moved to any edge of $t_p$. Specifically, it is enough to
consider the case where $J=(J_1,J_2)$, $I'=(I,J_1)$ and $J'=J_2$
so that $I\cdot(J_1,J_2)=(I,J_1)\cdot J_2$ as in
Figure~\ref{tree-move-trees-fig}. (Here we are assuming that $W_J$
is not order--$0$ since if both $W_I$ and $W_J$ are order--$0$ there
is nothing to prove.) The proof is given by the maneuver
illustrated in Figure~\ref{tree-move1-OR-fig}: Use the Whitney
disk $W_J$ to guide a Whitney move on $W_{J_1}$. This eliminates
the intersections between $W_{J_1}$ and $W_{J_2}$ (as well as
eliminating $W_J$ and $p$) at the cost of creating a new
cancelling pair of intersections between $W_{J_1}$ and $W_I$. This
new cancelling pair can be paired by a Whitney disk $W_{(I,J_1)}$
having a single intersection point $p'$ with $W_{J_2}$. That this
achieves the desired effect on the punctured tree can be seen in
Figure~\ref{tree-move2-OR-fig} by referring to the signs and
orientations in Figure~\ref{tree-move1-OR-fig}. See also the
discussion in pages 20--22 of \cite{ST} which includes group
elements.
\end{proof}
\newpage

\subsection{Algebraically- and geometrically-cancelling pairs}
\label{alg-geo-pairs}
Let $\widehat{\cT}_n(\pi,m)$ denote the group of order--$n$ decorated trees
modulo all the relations in Figure~\ref{Relations-fig} {\em except}
the IHX relation. We say that a pair of intersection points $p_+$
and $p_-$ in $\cW$ {\em cancel algebraically} if $\epsilon _{p_+}\cdot
t_{p_+}=-\epsilon _{p_-}\cdot t_{p_-}\in\widehat{\cT}_n(\pi,m)$. There
is a summation map that sends the disjoint union $t_n(\cW)=\amalg_p
\e_p\cdot t_p$ to an element $\widehat{\tau}_n(\cW):=\sum_p \e_p\cdot
t_p\in\widehat{\cT}_n(\pi,m)$ and the vanishing of $\widehat{\tau}_n(\cW)$
is equivalent to being able to arrange all of the order--$n$ intersection
points of $\cW$ into algebraically-cancelling pairs.

Given an algebraically-cancelling pair $p_{\pm}$ in a split Whitney
tower, one can chose orientations and whiskers on the Whitney disks in
the split subtowers containing $p_{\pm}$ so that the trees $t_{p_{\pm}}$
have identical orientations (and decorations) with $\e_{p_+}=-\e_{p_-}$.
(This is because the OR, HOL and AS relations are realized by these
choices, as described in Sections \ref{decorating-tree-subsec} and
\ref{subsec:AS}.)

A pair of intersection points $p_+$ and $p_-$ in $\cW$ {\em cancel
geometrically} if they can be paired by a Whitney disk. Geometric
cancellation implies algebraic cancellation, but the converse is not
true since two algebraically-cancelling intersection points might not
lie on the same Whitney disks.

The next lemma gives sufficient conditions for a sort of converse
involving some additional work.

\subsection{Simple intersection points and the transfer lemma}
\label{simple-int-transfer-lem-subsec}

\begin{figure}[ht!]
\centerline{\includegraphics[scale=.6]{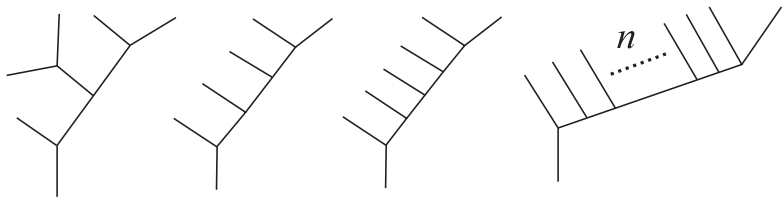}}
\caption{From left to right, the non-simple tree of lowest-order
(order--4) and the simple trees of order 4, 5
and $6+n$.}
\label{simple-trees-fig}
\end{figure}

Following the terminology of \cite{MKS} for iterated commutators of group
elements, we say that an intersection point $p\in \cW$ is {\em simple}
if its tree $t_p$ is simple (right- or left-normed) as illustrated in
Figure~\ref{simple-trees-fig}. The proof of the next lemma shows how to
exchange simple algebraically-cancelling pairs of intersection points
for geometrically-cancelling pairs.

\begin{lem}\label{transfer-lemma}
Let $\cW$ be an order--$n$ Whitney tower on order--0 surfaces $A_i$ such
that all order--$n$ intersection points of $\cW$ come in {\em simple}
algebraically-cancelling pairs. Then the $A_i$ are homotopic (rel
boundary) to $A'_i$ which admit an order--$(n+1)$ Whitney tower.
\end{lem}

\begin{figure}[ht!]
\centerline{\includegraphics[scale=.5]{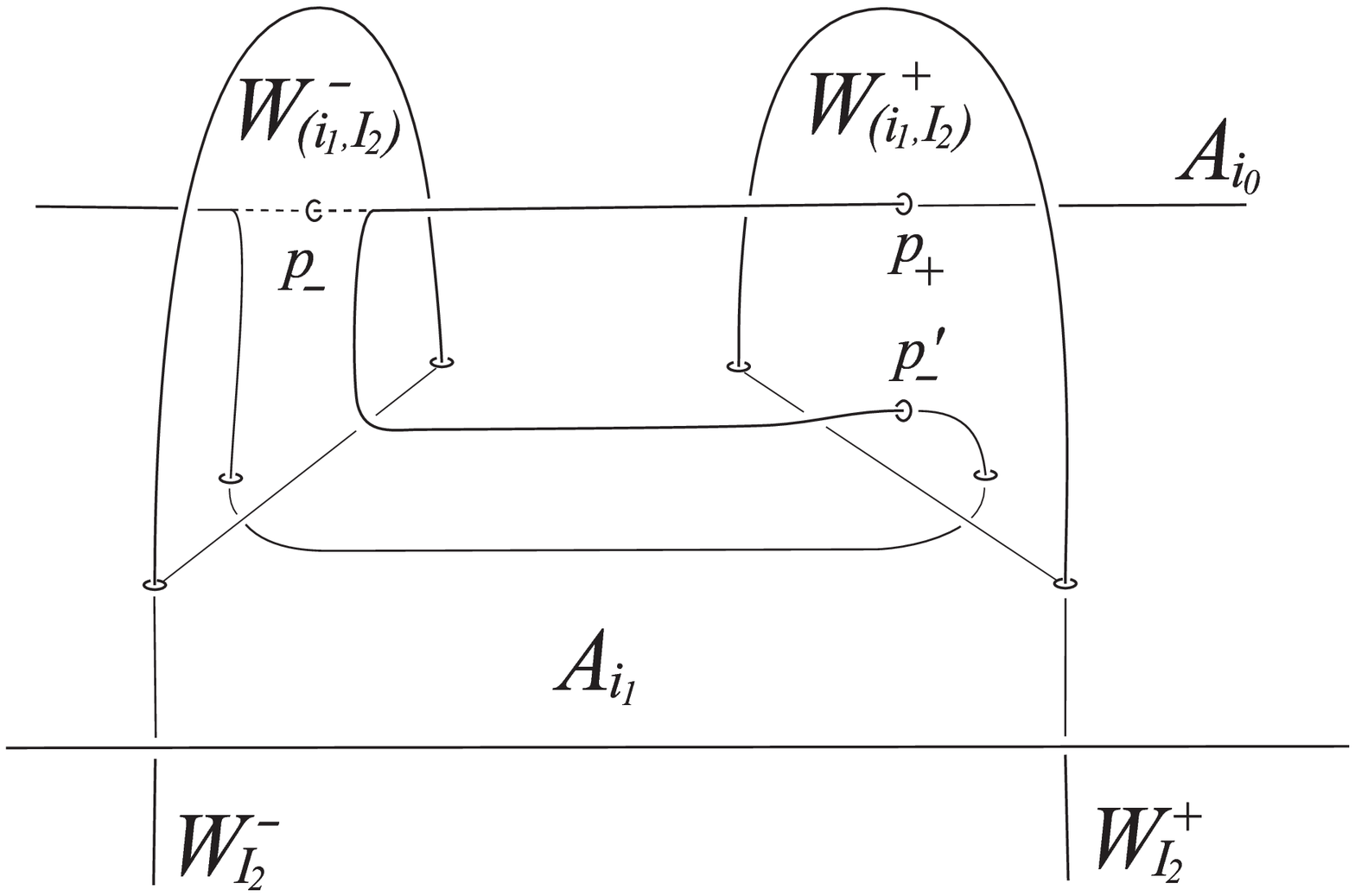}}
\nocolon \caption{}
\label{transfer-move-1A-fig}
\end{figure}

\begin{proof}
We will describe a modification of $\cW$ which exchanges one
algebraically-cancelling simple pair of order $n$ for another at the
cost of only creating geometrically-cancelling pairs. Iterating this
modification will, at the $n$th iteration, exchange an
algebraically-cancelling pair for {\em only} geometrically-cancelling
pairs. This modification is described in \cite{Y} for the case $n=1$
in a simply-connected manifold. (See also \cite{ST} for the $n=1$
non-simply-connected case.)  Applying this procedure to all
algebraically-cancelling pairs will complete the proof. We will discuss
only the simply-connected case; the reader can easily add group elements
to the figures (as in \cite{ST}).

We may assume that $\cW$ is split by Lemma~\ref{split-tower-lem}.
Let $p_+$ and $p_-$ be a simple algebraically-cancelling pair of
order--$n$ intersection points in $\cW$. By ``pushing the puncture out
to an end of the simple tree'' using Lemma~\ref{subtower-lemma},
we may further assume that $p_+$ and $p_-$ are intersections
between some order--0 surface $A_{i_0}$ and order--$n$ Whitney
disks $W^+_{I_1}$ and $W^-_{I_1}$ respectively where, for
this proof only, $I_k$ will denote a simple bracket of the form
$I_k:=(i_k,(i_{(k+1)},(\ldots,(i_{n},i_{(n+1)})\dots))\!=\!(i_k,I_{(k+1)})$
for $1\leq k\leq n+1$ and $I_{(n+1)}=i_{(n+1)}$.
\begin{figure}[ht!]
        \centerline{\includegraphics[scale=.5]{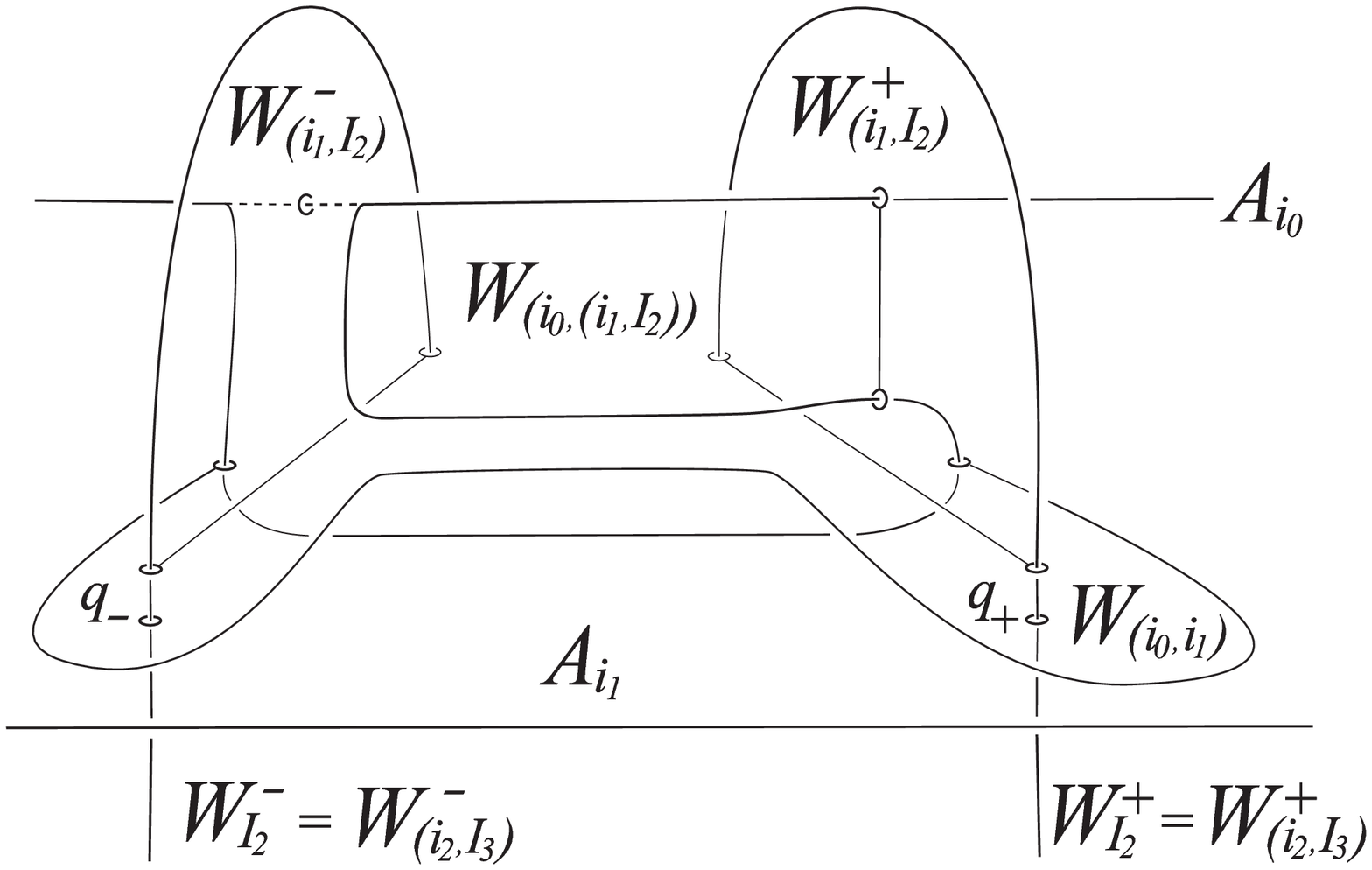}}
        \nocolon \caption{}
        \label{transfer-move-2A-fig}
\end{figure}

The first step in the modification is illustrated in
Figure~\ref{transfer-move-1A-fig} which shows how to exchange
$p_-\in A_{i_0}\cap W^-_{I_1}$ for $p'_-\in A_{i_0}\cap
W^+_{I_1}$, which cancels geometrically with $p_+$, at the cost of
creating a geometrically-cancelling pair of intersection points
between $A_{i_0}$ and $A_{i_1}$. Note that this first step is
possible because both $A_{i_0}$ and $A_{i_1}$ are {\em connected}.
The modification is completed by choosing Whitney disks for the
new geometrically-cancelling pairs as illustrated in
Figure~\ref{transfer-move-2A-fig}, which shows that a new
algebraically-cancelling pair $q_{\pm}\in W_{(i_0,i_1)}\cap
W^{\pm}_{I_2}$ has been created (recall that boundaries of Whitney
disks must be disjointly embedded). In the case $n=1$, $q_{\pm}$
would also cancel geometrically since then $I_{(n+1)}=i_{(n+1)}$
means that $W^{+}_{I_2}=W^{-}_{I_2}=A_{i_2}$ which is connected.
Note that $W_{(i_0,i_1)}$ is embedded (in a neighborhood of a
contractible 1--complex) and contains only the pair $q_{\pm}$ in
its interior. The Whitney disk $W_{(i_0,(i_1,I_2))}$ may intersect
anything but we don't care because it is a Whitney disk of order
$n+1$ and hence can only contain intersections of order strictly
greater than $n$.
\begin{figure}[ht!]
\centerline{\includegraphics[scale=.5]{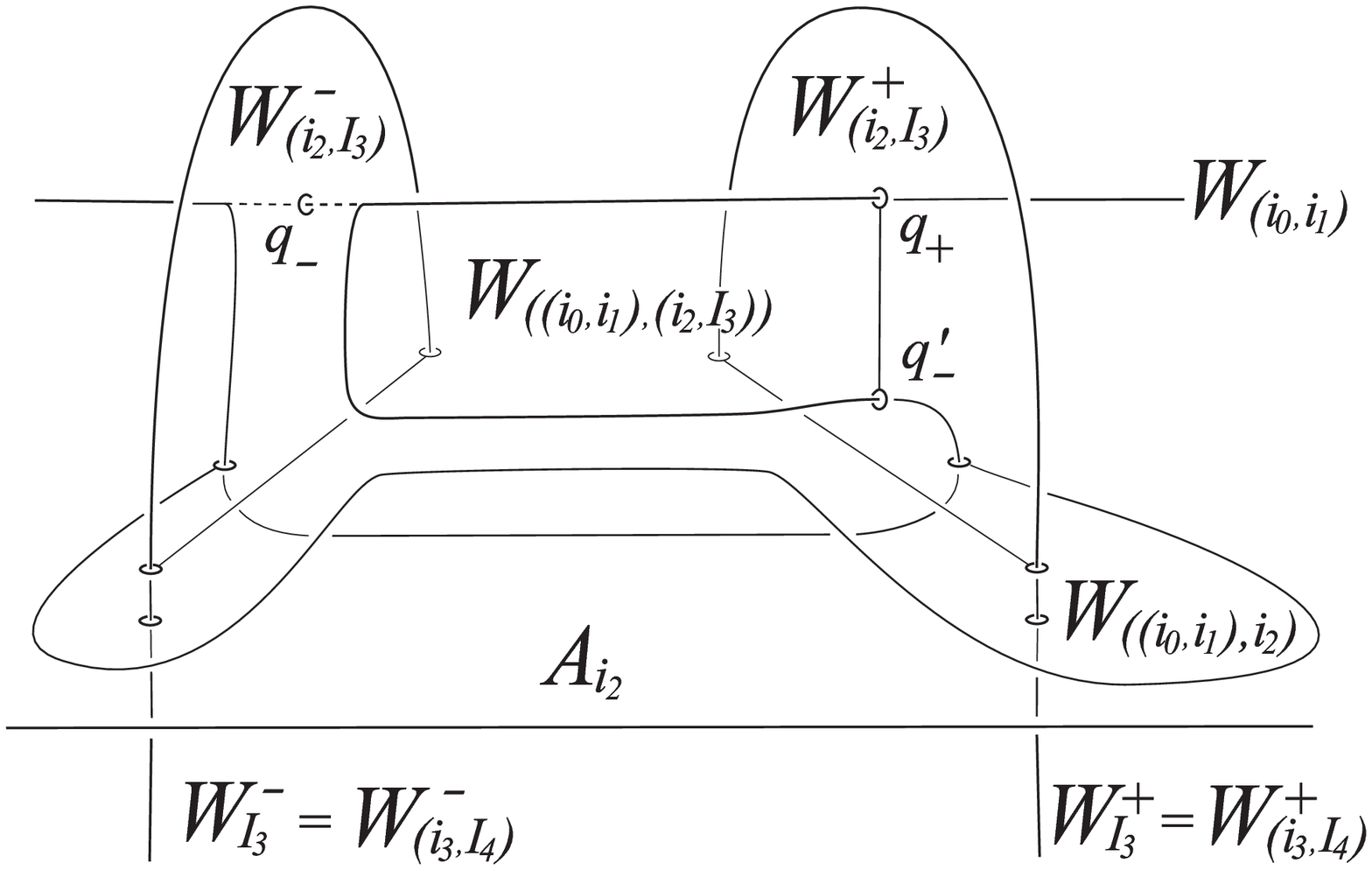}}
\nocolon \caption{}
\label{transfer-move-2B-fig}
\end{figure}
Now, assuming $n\geq 2$, apply this modification to $q_{\pm}$ as
illustrated in Figure~\ref{transfer-move-2B-fig}. Note that this
is only possible because we have the {\em connected} surface
$A_{i_2}$ to ``push along'', since we originally started with the
{\em simple} pair $p_{\pm}$ so that $W^{\pm}_{I_2}=W_{(i_2,I_3)}$.
\begin{figure}[ht!]
\centerline{\includegraphics[scale=.5]{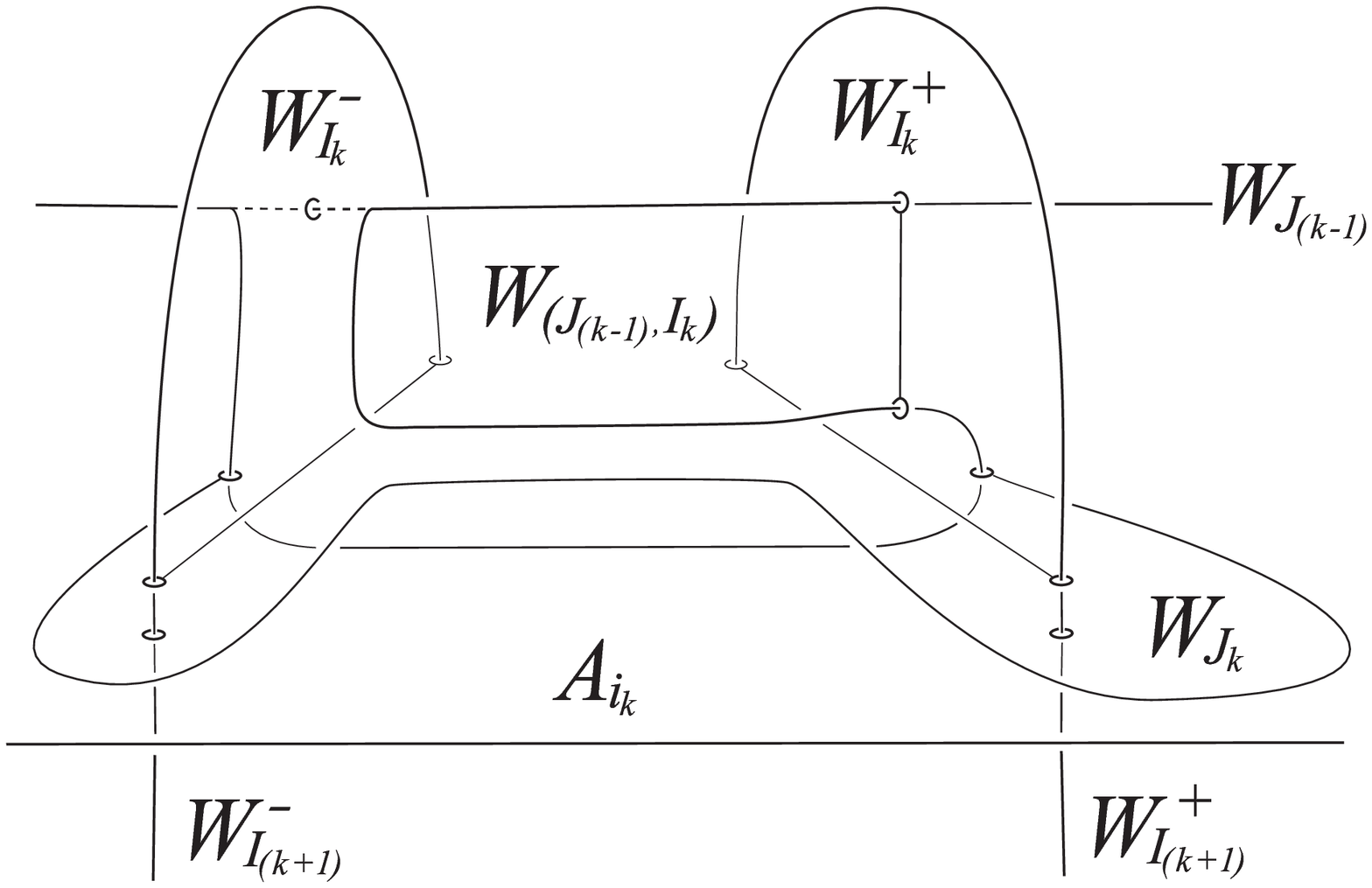}}
\nocolon \caption{}
\label{transfer-move-2C-fig}
\end{figure}
The $k$th iteration of this modification is illustrated in
Figure~\ref{transfer-move-2C-fig} where, for this proof only, we
denote the simple bracket $J_k:=(\ldots((i_0,i_1),i_2),\ldots,i_k)$
for $1\leq k\leq n$. The procedure terminates when $k=n$ meaning that
$W^{\pm}_{I_{(k+1)}}=W_{I_{(n+1)}}=A_{i_{(n+1)}}$ which is connected so
only geometrically-cancelling pairs are created.

This procedure can be applied to all the (simple) algebraically-cancelling
pairs: One can always find disjoint arcs between Whitney arcs
in the $A_{i_k}$ to guide the modification and all new Whitney disks of
order~$\leq n$ are contained in neighborhoods of these arcs so that no
unexpected intersections of order less than or equal to $n$ are created.
\end{proof}

\subsection{Geometric IHX and the Proof of Theorem~\ref{thm:build-tower}}
\label{build-thm-proof}
Given $\cW$ as in Theorem~\ref{thm:build-tower}, we will reduce the proof
to the case handled by Lemma~\ref{transfer-lemma} by using geometric
constructions and results from \cite{CST} and \cite{S}. Achieving the
hypotheses of Lemma~\ref{transfer-lemma} will involve two steps: First
$\cW$ will be modified to have only algebraically-cancelling pairs by
using the ``$4$--dimensional IHX construction'' in \cite{CST}. Then the
algebraically-cancelling pairs will be exchanged for simple
algebraically-cancelling pairs, using a related IHX construction of
\cite{S}. This second step is based on the effect of doing a Whitney move
on a Whitney disk in a split subtower and mimics the usual algebraic proof
that the group of unitrivalent trees modulo the IHX and AS relations is
spanned by simple trees (\cite{B2}, \cite{CT1}).

\subsection{Creating algebraically-cancelling pairs}
The vanishing of $\tau_n(\cW)\in{\cT}^t_n(\pi,m)$
means that $\tau_n(\cW)$ lifts to
$\widehat{\tau}_n(\cW)\in\mathrm{span}\{\mathrm{I}-\mathrm{H}+\mathrm{X}\}
<\widehat{\cT}_n(\pi,m)$.
To get only algebraically-cancelling pairs we apply the following
corollary of the {\em $4$--dimensional IHX Theorem} in \cite{CST}:
\begin{prop}\label{prop:IHX} Let $\cW$ be any order--$n$ Whitney tower on
order--0 surfaces $A_i$. Then, given any decorated order--$n$ unitrivalent
trees $\mathrm{I}$, $\mathrm{H}$ and $\mathrm{X}$ differing only by the
local $\mathrm{IHX}$ relation of Figure~\ref{Relations-fig}, there exists
an order--$n$ Whitney tower $\cW'$ on $A'_i$ homotopic (rel boundary)
to the $A_i$ such that
$$t_n(\cW')=t_n(\cW)+\mathrm{I}-\mathrm{H}+\mathrm{X}.$$
\end{prop}
Note that the ``sum'' on the right hand side is really a disjoint
union of signed decorated trees; the summation map takes this
equation to the corresponding equation in
$\widehat{\cT}_n(\pi,m)$.

\begin{proof}
As observed in Remark~\ref{rem:w-tree-f-move}, creating a ``clean''
Whitney disk by applying a finger move to surfaces in a Whitney tower
``realizes'' the rooted product $\ast$ on the corresponding rooted
trees. Since finger moves are supported near arcs, one can modify $\cW$
to create any number of clean Whitney disks realizing arbitrary rooted
decorated trees without changing $t_n(\cW)$. Let $W^i$, $i=1,2,3,4$ be
four such Whitney disks which correspond to the four fixed vertices of
the trees I, H and X in the statement. (Of course if any of the fixed
vertices is univalent then the corresponding ``Whitney disk'' is just
an order--0 surface.)

Now the $4$--dimensional IHX Theorem of \cite{CST} says that there
exists an order--2 Whitney tower $\cW_{\mathrm{IHX}}$ on oriented
$2$--spheres $A_i$, $i=1,2,3,4$, in a $4$--ball having geometric
intersection tree $t_2(\cW_{\mathrm{IHX}})$ equal precisely to the
order--2 IHX relation. So by tubing $A_i$ into $W^i$, for each $i$,
we can get $\cW'$ as desired. No unexpected intersections are
created since the entire construction takes place near a
collection of arcs and the (arbitrarily small) $4$--ball. (In the
decorated case the desired group elements are controlled by the
tubes.)
\end{proof}

So by applying Proposition~\ref{prop:IHX} as necessary we can assume
that $\widehat{\tau}_n(\cW)=0\in\widehat{\cT}_n(\pi,m)$ which means
that all order--$n$ intersection points can be arranged in
algebraically-cancelling pairs.

\subsection{Simplifying the cancelling pairs}
In case there are algebraically-cancelling pairs which are not simple, we
appeal to results in \cite{S}: Proposition~7.1 of \cite{S} describes an
algorithm for modifying a Whitney tower to have only simple intersection
points. This geometric algorithm, which mimics the algebraic algorithm
described in \cite{B2} and \cite{CT1}, depends on a ``Whitney move''
version of the IHX relation (Lemma~7.2 of \cite{S}) which replaces
a split subtower $\cW_p$ by two split subtowers $\cW_{p'}$ and
$\cW_{p''}$ and has the effect of replacing $\e_p\cdot t_p=\mathrm{I}$
by $\mathrm{H}-\mathrm{X}=\e_{p'}\cdot t_{p'}+\e_{p''}\cdot t_{p''}$
in the geometric intersection tree.  The point of the algorithm is
that the trees H and X are ``closer'' to being simple and by iterating
one is eventually left with only simple trees. (The construction is
supported in a neighborhood of $\cW_p$ so no unwanted intersections are
created.) Although Proposition~7.1 and Lemma~7.2 of \cite{S} are only
proved in the unoriented undecorated case it is not hard to add signs
to the intersection points in the diagrams in \cite{S} and apply the
conventions of this paper, especially having seen the related proof of
Lemma~\ref{transfer-lemma} above.

So in the present setting we have only algebraically-cancelling
pairs of order--$n$ intersection points in an order--$n$ Whitney
tower $\cW$ which we may assume is split by
Lemma~\ref{split-tower-lem}. If any of these cancelling pairs are
not simple, then we apply the just-mentioned IHX algorithm of
\cite{S} {\em pairwise} (so as to preserve
$\widehat{\tau}_n(\cW)=0\in\widehat{\cT}_n(\pi,m)$) until we are
left with only simple algebraically-cancelling pairs. The proof of
Theorem~\ref{thm:build-tower} is now complete by
Lemma~\ref{transfer-lemma}.\hfill$\square$


\section{Proof of Theorem~\ref{thm:Milnor}}

The proof of Theorem~\ref{thm:Milnor} uses results from \cite{HM},
\cite{KT} and \cite{S} as well Theorem~\ref{thm:build-tower} to compare
an arbitrary link $L$ to certain well-known standard links which generate
the first non-vanishing Milnor and $Z^t$ invariants.

\begin{figure}[ht!]
\centerline{\includegraphics[scale=.35]{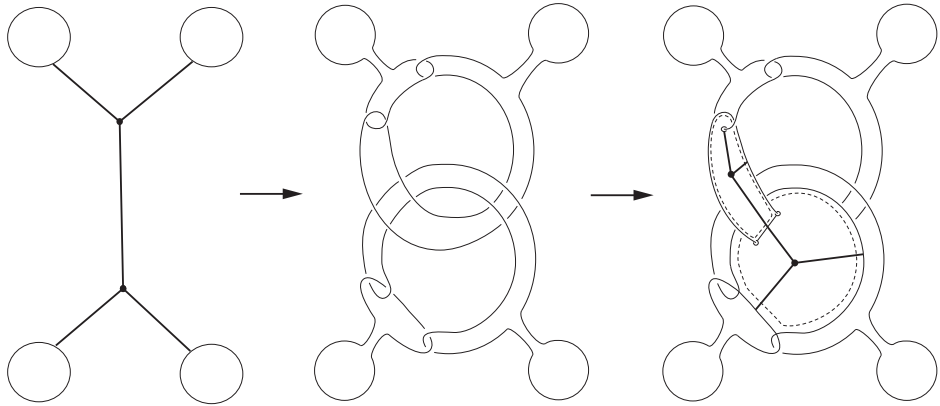}}
\caption{From left to right: An order--2 (positively signed)
vertex-oriented tree I whose univalent vertices correspond to the
components of an unlink; The Bing--Cochran--Habiro link $L_\mathrm{I}$;
An order--2 Whitney tower $\cW$ bounded by $L_\mathrm{I}$ with
$\tau_2(\cW)=\mathrm{I}$.}
\label{Cochran-Habiro-fig}
\end{figure}

\subsection{Bing--Cochran--Habiro links}\label{Cochran-Habiro-subsec}
Given a collection $\sigma$ of signed labelled vertex-oriented
order--$n$ trees, Cochran \cite{C} and Habiro \cite{H} have described,
using Bing doubling and clasper surgery respectively, how to construct
(from the unlink) a link $L_{\sigma}$ such that $K_n(L_{\sigma})$ is
represented by $\sigma$ (considered as a sum). Habiro's construction
applies more generally to unitrivalent {\em graphs}, but for trees the
two constructions coincide (by applying Kirby calculus to a framed link
surgery description).

Given such a {\em Bing--Cochran--Habiro link} $L_{\sigma}$, we will use the following two facts:
\begin{enumerate}
\item $L_{\sigma}$ bounds an order--$n$ Whitney tower $\cW_{\sigma}$
with $\tau_n(\cW_{\sigma})=\sigma\in \cT_n(m)$.
\item $K_n(L_{\sigma})=\sigma\in \cT_n(m)\otimes\Q$.
\end{enumerate}
The Whitney tower $\cW$ in statement (i) is easily constructed by
``pulling apart'' a Bing double in Cochran's construction (see
Figure~\ref{Cochran-Habiro-fig}): This creates Whitney disks whose
boundaries are essentially the {\em derived links} in \cite{C} and each
$t_p\in\sigma$ corresponds to a {\em derived linking}. Alternatively,
starting with Habiro's clasper surgery description one can apply
the translation to {\em grope cobordism} of \cite{CT1} and then the
translation to Whitney towers of \cite{S} and \cite{CST}.

For statement (ii), see Section 8 of \cite{HM}. Although \cite{HM} works
with {\em string} links, the first non-vanishing term of $Z^t(L)-1$
is equal to the first non-vanishing term of $Z^t(SL)-1$ where $SL$
is any string link whose closure is $L$ (see Section~5 of \cite{MV}).

\subsection{Whitney towers and the Kontsevich integral}
Let $L$ and $\cW$ be as hypothesized in Theorem~\ref{thm:Milnor}.
Denote by $\sigma$ any disjoint union of signed (labelled
vertex-oriented) trees which represents $\tau_n(\cW)\in\cT_n(m)$,
eg the geometric intersection tree $t(\cW)$ of $\cW$
(\ref{geo-int-tree-w-tower}). Let $L_{\sigma}$ be a
Bing--Cochran--Habiro link formed from the unlink using $\sigma$.
Then, by (i) of \ref{Cochran-Habiro-subsec}, $L_{\sigma}$ bounds
an order--$n$ Whitney tower $\cW_{\sigma}$ in $D^4$ with
$\tau_n(\cW_{\sigma})=\tau_n(\cW)\in \cT_n(m)$. Now think of $\cW$
and $\cW_{\sigma}$ as each sitting in a copy of $S^3\times I$
($D^4$ with a neighborhood of a point removed). By gluing together
the two copies of $S^3\times I$ (along the $S^3$ boundary of the
removed neighborhoods) and connecting each order--0 $2$--disk of
$\cW$ with the corresponding order--0 $2$--disk of
$\cW_{\sigma}$ by a small tube we get properly immersed annuli
$A_i$ in $S^3\times I$ cobounded by the link components. Since the
tubes may be chosen to avoid creating new intersection points, the
$A_i$ admit an order--$n$ Whitney tower $\cW'$ with
$$\tau_n(\cW')=\tau_n(\cW)-\tau_n(\cW_{\sigma})=0\in \cT_n(m)$$
where the minus sign comes from reversing the orientation of one of
the two copies of $S^3\times I$. By Theorem~\ref{thm:build-tower},
the vanishing of $\tau_n(\cW')$ implies that (after a homotopy rel
boundary) the $A_i$ admit a Whitney tower of order~$n$, that is, $L$ and
$L_{\sigma}$ are {\em order--$n$ Whitney equivalent}. By the main theorem
in \cite{S}, order--$n$ Whitney equivalence implies (in fact is equivalent
to) {\em class~$(n+1)$ grope concordance}, meaning that we can conclude
that the components of $L$ and $L_{\sigma}$ cobound disjoint properly
embedded annulus-like gropes of class~$(n+1)$. This implies, by \cite{KT}
Corollary~4.2, that $L$ and $L_{\sigma}$ have the same $\mu$--invariants
of length less than or equal to $(n+1)$. It follows from \cite{HM} that
$K_n(L)=K_n(L_{\sigma})$ which is equal to $\sigma\in\cT_n(m)\otimes\Q$
by (ii) of \ref{Cochran-Habiro-subsec} above.
\endproof

\Addresses\recd

\end{document}